%% file: CMofBessel_article.tex
\documentclass[11pt]{article}
\usepackage[english]{babel}
\usepackage{amssymb}
\usepackage{amsmath}
\usepackage{amsthm}
\usepackage[a4paper,margin=2.54cm]{geometry}
\usepackage{graphicx}
\usepackage{todonotes}
\usepackage{comment}
\usepackage{appendix}
\usepackage{tikz}
\usetikzlibrary{decorations.pathreplacing}
\usetikzlibrary{matrix}
\usepackage{pgfopts}
\usepackage{array}
\usepackage{mathtools}
\usepackage{ytableau}
\usepackage{tikz-cd}

\theoremstyle{plain}
\newtheorem*{theorem*}{Theorem}
\newtheorem{theorem}{Theorem}[section] 
\newtheorem{lemma}[theorem]{Lemma}
\newtheorem{proposition}[theorem]{Proposition}
\newtheorem{corollary}[theorem]{Corollary}

\theoremstyle{definition}

\newtheorem{remark}[theorem]{Remark}

\makeatletter
\@addtoreset{equation}{section}
\makeatother
\numberwithin{equation}{section}

\def\NEARROWK{{%
		\setbox0\hbox{$\nearrow$}%
		\rlap{\hbox to \wd0{\hss $\phantom{x}_k$}}\box0
}}

\def\NEARROWKTILDE{{%
		\setbox0\hbox{$\nearrow$}%
		\rlap{\hbox to \wd0{\hss $_{\tilde{k}}$}}\box0
}}

\def\NEARROWKtwo{{%
		\setbox0\hbox{$\nearrow$}%
		\rlap{\hbox to \wd0{\hss $\phantom{x}_{2}$}}\box0
}}

\def\SEARROWK{{%
		\setbox0\hbox{$\searrow$}%
		\rlap{\hbox to \wd0{\hss $\phantom{x}^k$}}\box0
}}

\usepackage{authblk}
\title{Universality for conditional measures of the Bessel point process}

\author[1]{Leslie Molag}
\author[2]{Marco Stevens}
\affil[1,2]{KU Leuven, Department of Mathematics, 
	
	E-mail: {\tt leslie.molag@kuleuven.be} and {\tt marco.stevens@kuleuven.be}
}

\date{}                     
\providecommand{\keywords}[1]{\textit{\textit{Keywords: }} #1}

\usepackage{hyperref}
\hypersetup{colorlinks, citecolor=black, linkcolor=black}
\AtBeginDocument{} 

\begin{document}
\maketitle

\input{abstract}
\keywords{Bessel point process, rigidity, conditional measures, orthogonal polynomial ensembles, asymptotics, Riemann-Hilbert analysis.
}

\input{mainresult}

\input{motivation}
\input{behaviour}

\input{outlineofproof}

\input{shortRHproblem}

\begin{appendices}

\end{appendices}

\section*{Acknowledgements}

Leslie Molag is supported by a PhD fellowship of the Flemish Science Foundation (FWO). Marco Stevens is supported by EOS project 30889451 of the Flemish Science Foundation (FWO) and by the Belgian Interuniversity Attraction Pole P07/18. Both authors are partly supported by the long term structural funding-Methusalem grant of the Flemish Government. They are also grateful to Arno Kuijlaars for useful discussions and for proofreading the article.

\end{document}

%% file: abstract.tex
\begin{abstract}\label{abstract}
The Bessel point process is a rigid point process on the positive real line and its conditional measure on a bounded interval $[0,R]$ is almost surely an orthogonal polynomial ensemble. In this article, we show that if $R$ tends to infinity, one almost surely recovers the Bessel point process. In fact, we show this convergence for a deterministic class of probability measures, to which the conditional measure of the Bessel point process almost surely belongs.
\end{abstract}

%% file: mainresult.tex
\section{Main result}
\label{sec:mainresult}
In this article, we are concerned with sequences $0 < p_1 < p_2 < \cdots$ that satisfy the growth condition
\begin{equation}
\label{eq:growthconditionpoints}
\lim_{n\rightarrow \infty} \frac{p_n}{n^2}=\pi^2.
\end{equation}
Our major motivating example of such sequences is formed by the \textbf{Bessel point process}; this is the determinantal point process whose kernel is the \textbf{Bessel kernel}
\begin{equation}
\label{eq:Besselkernel}
J_\nu(x,y)=\frac{J_\nu(\sqrt{x})\sqrt{y}J'_\nu(\sqrt{y}) - J_\nu(\sqrt{y})\sqrt{x}J'_\nu(\sqrt{x})}{2(x-y)}, \qquad x,y \in (0,\infty),
\end{equation} 
where $\nu>-1$ is a parameter and we (ab)use the notation $J_\nu$ to mean the Bessel function if it has one variable and the Bessel kernel if it has two variables. The configurations taken from the Bessel point process almost surely satisfy the growth condition \eqref{eq:growthconditionpoints}; see Proposition \ref{prop:behaviourBessel}. 

The aim of this paper is to study the asymptotic behaviour of the conditional measure \cite{Bufetov_conditionalmeasures,Kuijlaars_MinaDiaz} of the Bessel point process; see Section \ref{sec:motivation} for more details. To this extent, we define weights $\bar{w}_{X,\nu,R}$ on the interval $[0,R]$, where $X=(p_n)_{n=1}^\infty$ represents an increasing sequence of positive numbers satisfying the growth condition \eqref{eq:growthconditionpoints}, $R>0$ and $\nu>-1$. This weight is defined by
\begin{equation}
\label{eq:weight}
\bar{w}_{X,\nu,R}(t)=t^\nu \prod_{p_n>R} \Big(1-\frac{t}{p_n}\Big)^2, \qquad t\in [0,R].
\end{equation}
We study the finite point process associated to these weights. More specifically, we denote, for every $R>0$, by $N(R)$ the number of points in $X \cap [0,R]$ and we consider the point process with $N(R)$ points $(t_1,\dots,t_{N(R)})$ that all lie in $[0,R]$, with joint probability density function
\begin{equation}
\label{eq:jointpdf}
\frac{1}{Z} \prod_{1\leq i < j\leq N(R)} (t_i - t_j)^2 \prod_{i=1}^{N(R)} \bar{w}_{X,\nu,R}(t_i),
\end{equation}
where $Z$ is some normalization constant that depends on $X$, $\nu$ and $R$. If $X$ is taken as a configuration from the Bessel point process with parameter $\nu$, then the conditional measure of the Bessel point process associated to $X$ \cite{Bufetov_conditionalmeasures} is almost surely the point process described by \eqref{eq:jointpdf}. 

We are interested in the asymptotics of this point process as $R\rightarrow \infty$. For this, convergence is understood in the sense of \textit{convergence of kernels}. It is well-known that a point process with joint probability density function of the form \eqref{eq:jointpdf} is an orthogonal polynomial ensemble. In other words, it is a determinantal point process whose kernel is the normalized Christoffel-Darboux kernel (see Section \ref{sec:outline}) associated to the weight $\bar{w}_{X,\nu,R}$. For these kernels, we have the following asymptotic behaviour.

\begin{theorem}[Main Theorem]
\label{thm:maintheorem}
Suppose that $X=(p_n)_{n=1}^\infty$ is a strictly increasing sequence of positive numbers satisfying the growth condition \eqref{eq:growthconditionpoints}, let $\nu>-1$ and write $N(R)$ for the number of points in $X\cap [0,R]$ for all $R>0$. Then the normalized Christoffel-Darboux kernel associated to the weights $\bar{w}_{X,\nu,R}$ (as defined in \eqref{eq:weight}) satisfies
\begin{equation}\label{eq:limitkernel}
\lim_{R\to \infty} K_{N(R)}(x,y; \bar{w}_{X,\nu,R}) = J_\nu(x,y), \qquad (x,y)\in (0,\infty)^2
\end{equation}
uniformly on compact sets, where $J_\nu$ is the Bessel kernel with parameter $\nu$, as given in \eqref{eq:Besselkernel}.
\end{theorem}
Since this limit is independent of the exact choice of $X$, we speak of \textbf{universality} of the conditional measure. In this way, this article can be seen as a continuation of the result for the conditional measures of the \textit{sine point process} obtained in \cite{Kuijlaars_MinaDiaz}.

\begin{remark}
We note that if one alters the growth condition \eqref{eq:growthconditionpoints} to $\lim_{n\rightarrow \infty} \frac{p_n}{n^2}=c\pi^2$ for some $c>0$ that one obtains the rescaled Bessel kernel $\frac{1}{c} J_\nu\left(\frac{x}{c},\frac{y}{c}\right)$ as the limit in \eqref{eq:limitkernel}. Hence the real importance of the growth condition \eqref{eq:growthconditionpoints} is the quadratic nature of $p_n$; the constant $\pi^2$ can be seen as an artifact of the conventional scaling choice of the Bessel functions and kernels.
\end{remark}

The rest of this article is dedicated to proving Theorem \ref{thm:maintheorem} and is outlined as follows:
\begin{itemize}
\item In Section \ref{sec:motivation}, we give the reason of interest for this type of result. We discuss the notions of rigidity and conditional measures, which should be seen as the enveloping framework for this article.
\item In Section \ref{sec:behaviour}, we show that the configurations taken from the Bessel point process almost surely satisfy the growth condition \eqref{eq:growthconditionpoints}.
\item In Section \ref{sec:outline}, we give the proof of Theorem \ref{thm:maintheorem}, except for a Riemann-Hilbert analysis, which we give in Section \ref{sec:RHproblem}. The structure of the proof is analogous to the one given in \cite{Kuijlaars_MinaDiaz}, which solves the similar question for the conditional measure of the sine point process.
\end{itemize} 

%

%% file: motivation.tex
\section{Motivation}
\label{sec:motivation}

This paper should be considered as contributing to the understanding of \textit{rigid} point processes. Surrounding this concept of rigidity, there are three subsequent questions for a given point process:
\begin{enumerate}
\item Is the point process rigid?
\item If the point process is rigid, what is the conditional measure of the point process on a compact subspace? Can one give closed formulas for these conditional measures for a family of well-chosen compact subsets?
\item If one has the answer to the previous question, what is the asymptotic limit of the conditional measure as the compact subset grows to cover the whole space? 
\end{enumerate}
These questions have been studied in the literature for several (classes of) point processes. See for example \cite{Bufetov_rigidity, Ghosh_Peres} for the first, \cite{Bufetov_conditionalmeasures} for the second and \cite{Kuijlaars_MinaDiaz} for the third question. We refer the reader to these sources for an extensive discussion about what rigidity and conditional measures are; in Section \ref{sec:mainresult} we stated everything we need in this paper from these concepts.

In this article, we answer the above question 3 for the Bessel point process, i.e. the point process whose kernel is the Bessel kernel \eqref{eq:Besselkernel}. For this process, the first question above has been answered positively in \cite{Bufetov_rigidity}. Subsequently, in \cite{Bufetov_conditionalmeasures}, the second question has been answered explicitly for compact subsets of the form $[0,R]$. As mentioned in Section \ref{sec:mainresult}, the conditional measure of the Bessel point process almost surely is an orthogonal polynomial ensemble, with \eqref{eq:weight} as the involved weight. 

The natural next question for the Bessel point process is then the third question above. It is clear that one needs to understand the behaviour of the points $(p_n)_{n=1}^\infty$ as $n\rightarrow \infty$ in a configuration taken from the Bessel point process in order to be able to understand the asymptotics of the weight \eqref{eq:weight} as $R\rightarrow \infty$. For this, we have the following result.
\begin{proposition}
\label{prop:behaviourBessel}
Suppose that the points $0 < p_1 < p_2 < \cdots$ form a configuration taken from the Bessel point process with some parameter $\nu>-1$. Then for every $\epsilon>0$ we have, with probability $1$, that
\begin{equation}
\label{eq:behaviourBessel}
p_n = \pi^2 n^2 + \mathcal O\left(n \sqrt n \log^{1+\epsilon} n\right), \qquad \text{as }n\to\infty.
\end{equation}
\end{proposition}
We prove this result in Section \ref{sec:behaviour}. By this Proposition, we know that the configurations taken from the Bessel point process almost surely satisfy the growth condition \eqref{eq:growthconditionpoints}. Hence our main result Theorem \ref{thm:maintheorem} answers the above third question for the Bessel point process; as one would expect, the asymptotic limit of the conditional measure on the compact subset $[0,R]$ is the Bessel point process itself, almost surely with respect to the chosen configuration $X$. In fact, Theorem \ref{thm:maintheorem} shows that these asymptotics are even more universal, since the growth condition \eqref{eq:growthconditionpoints} does not require any big-O term like the one appearing in \eqref{eq:behaviourBessel}.

%% file: behaviour.tex
\section{Proof of Proposition \ref{prop:behaviourBessel}}
\label{sec:behaviour}

For the proof of Proposition \ref{prop:behaviourBessel}, we rely on known results about the Bessel process. For this, we let $(p_n)_{n=1}^\infty$ be a configuration taken from the Bessel point process such that $0 < p_1 < p_2 <\cdots$, and for every $T>0$, we let $N(T)$ be the random variable counting the number of points in $(p_n)_{n=1}^\infty \cap [0,T]$. Then, in \cite{Soshnikov}, it was shown that this variable has the following behaviour:
\begin{align}\label{eq:ETasymptotics}
\mathbb E N(T) &= \frac{\sqrt{T}}{\pi}+\mathcal O(1), & \qquad \mathrm{as} \ T\rightarrow \infty,\\ \label{eq:VarTasymptotics}
\operatorname{Var} N(T) &= \frac{1}{4\pi^2} \log(T)+\mathcal O(1), &\qquad \mathrm{as} \ T\rightarrow \infty.
\end{align}
We use this behaviour to prove Proposition \ref{prop:behaviourBessel}.
\begin{proof}[Proof of Proposition \ref{prop:behaviourBessel}]
It will be somewhat more practical to prove that as $n\to\infty$
\begin{equation}
\label{eq:behaviourpn+1}
p_{n+1}=\pi^2 n^2+\mathcal O\left(n\sqrt n \log^{1+\epsilon} n\right)
\end{equation}
with probability $1$, which is equivalent to the statement of the proposition. We set 
\[T_n = \pi^2 \left(n^2 + n \sqrt n \log^{1+\epsilon} n\right).\] For $n$ big enough we have
\begin{align} \label{eq:ETasymptotics1}
\frac{\sqrt{T_n}}{\pi} - n \geq \frac{1}{3} \sqrt n \log^{1+\epsilon} n\quad \text{ and } \quad \left|\frac{\sqrt{T_n}}{\pi} - \mathbb E N(T_n)\right| \leq \frac{1}{2}\left(\frac{\sqrt{T_n}}{\pi} - n\right),
\end{align}
where the first inequality directly follows from the definition of $T_n$ and the second inequality follows from the first and \eqref{eq:ETasymptotics}. Now let us suppose that $n$ is big enough in this sense, and that we have
\begin{align*}
\left|N(T_n) - \sqrt{T_n}/\pi\right| \geq \frac{\sqrt{T_n}}{\pi} - n.
\end{align*}
We can then invoke \eqref{eq:ETasymptotics1} and the reverse triangle inequality to obtain
\begin{align*}
\left|N(T_n) - \mathbb E N(T_n)\right| &\geq |N(T_n) - \sqrt{T_n}/\pi| - \left|\sqrt{T_n}/\pi - \mathbb E N(T_n)\right|\\
&\geq \frac{1}{2}\left(\frac{\sqrt{T_n}}{\pi} - n\right)
\geq \frac{1}{6} \sqrt n \log^{1+\epsilon} n.
\end{align*}
Next we make the observation that $p_{n+1} > T_n$ is equivalent to $N(T_n)\leq n$. Hence 
\begin{align*}
\mathbb P\left(p_{n+1}>T_n\right) &= \mathbb P\left(N(T_n) \leq n\right)\\
&\leq \mathbb P\left(\left|N(T_n) - \sqrt{T_n}/\pi\right|\geq \sqrt{T_n}/\pi - n\right)\\
&\leq \mathbb P\left(\left|N(T_n) - \mathbb E N(T_n)\right| \geq \frac{1}{6} \sqrt n \log^{1+\epsilon} n\right).
\end{align*}
Using Chebyshev's inequality we then obtain, for $n$ big enough
\begin{equation}
\label{eq:usingChebyshev}
\mathbb P\left(p_{n+1} > T_n\right)
\leq \frac{\operatorname{Var} N(T_n)}{\operatorname{Var} N(T_n) + \frac{1}{36} n \log^{2+2\epsilon} n}
\leq \frac{2}{n \log^{1+2\epsilon} n},
\end{equation}
where we have used \eqref{eq:VarTasymptotics} in the last step. We notice that
\[
\sum_{n=3}^\infty \frac{1}{n \log^{1+2\epsilon} n} \leq \int_2^\infty \frac{dt}{t \log^{1+2\epsilon} t} = \int_{\log 2}^\infty \frac{dt}{t^{1+2\epsilon}} < \infty,
\]
which means that the probability that appears in the left hand side of \eqref{eq:usingChebyshev} is summable, and hence the Borel-Cantelli lemma implies that, with probability $1$
\[
p_{n+1} > T_n = \pi^2 \left(n^2 + n \sqrt n \log^{1+\epsilon} n \right)
\]
occurs only finitely many times.  Analogous reasoning leads to a corresponding statement for a lower bound for $p_{n+1}$, and the proposition follows.
\end{proof}

%% file: outlineofproof.tex
\section{Proof of Theorem \ref{thm:maintheorem}}
\label{sec:outline}
In this section, we give the main components of the proof of Theorem \ref{thm:maintheorem}. The only element that we do not give in this section is the Riemann-Hilbert analysis needed to prove the asymptotics of the kernel associated to a weight defined in this section; we postpone this to the next section.

\subsection{Notation and general remarks}
We start with some conventions regarding notation. To any weight $w$ on a compact interval $[a,b]$, we denote the associated orthogonal polynomials by $(\varphi_j(\cdot; w))_{j=0}^\infty$, i.e., for all $i,j\geq 0$, we have
\begin{equation}
\label{eq:orthpol}
\int_{a}^b \varphi_j(t; w) \varphi_i(t; w) w(t) dt = \delta_{i,j}.
\end{equation}
These requirements do not completely define these polynomials. We use the convention that these polynomials have positive leading coefficients, which does define them uniquely. The \textbf{Christoffel-Darboux kernel} associated to the weight $w$ is defined by
\begin{equation}
\label{eq:defordCDkernel}
\widehat{K}_n(x,y; w) = \sum_{i=0}^{n-1} \varphi_i(x; w) \varphi_i(y; w), \qquad x,y\in [a,b].
\end{equation}
A slight modification leads to the \textbf{normalized Christoffel-Darboux kernel}
\begin{equation}
\label{eq:defnorCDkernel}
K_n(x,y; w) = \sqrt{w(x) w(y)} \sum_{i=0}^{n-1} \varphi_i(x; w) \varphi_i(y; w), \qquad x,y\in [a,b].
\end{equation}
It is well-known that the Christoffel-Darboux kernel has the property that
\begin{equation}
\label{eq:minpropChr}
\frac{1}{\widehat{K}_n(x,x; w)} = \min_{\deg P<n} \frac{1}{P(x)^2} \int_{a}^{b} \lvert P(t) \rvert^2 w(t) dt,
\end{equation}
where the quantity on the left hand side is called the \textit{Christoffel function} of the weight $w$. In the rest of this section, we will also need the simple behaviour of orthogonal polynomials (and the kernels built out of them) under rescalings. For easy reference, we summarize this behaviour in the following lemma.
\begin{lemma}
\label{lem:rescalingweights}
Suppose that $w$ is a weight on an interval $[a,b]$. Let $c,d>0$ and consider the new weight $\bar{w}(t):=dw(ct)$ on the interval $[a/c,b/c]$. Then we have that
\begin{equation}
\label{eq:orthpolrescaled}
\varphi_i(t;\bar{w}) = \frac{\sqrt{c}}{\sqrt{d}} \varphi_i(ct;w), \qquad t\in [a/c,b/c],
\end{equation}
for all integers $i\geq 0$. Furthermore, we have
\begin{equation}
\label{eq:norCDkernrescaled}
\frac{1}{c} K_n(\tfrac{x}{c},\tfrac{y}{c}; \bar{w}) = K_n(x,y;w), \qquad x,y\in [a,b],
\end{equation}
and
\begin{equation}
\label{eq:ordCDkernrescaled}
\frac{1}{c} \widehat{K}_n(\tfrac{x}{c},\tfrac{y}{c}; \bar{w}) = \frac{1}{d}\widehat{K}_n(x,y;w), \qquad x,y\in [a,b].
\end{equation}
\end{lemma}
One immediately checks this by rewriting the orthogonality relations.

\subsection{Rescaling to a point process on a fixed interval}
In order to prove Theorem \ref{thm:maintheorem}, we want to compare the considered point processes on $[0,R]$ for various $R$. For asymptotic analysis, it is more convenient to rescale these point processes to the interval $[0,1]$, such that all the considered point processes are defined on the same space. We define the new weight $w_{X,\nu,R}$ on $[0,1]$ by
\begin{equation}
\label{eq:weighton01}
w_{X,\nu,R}(t)=t^\nu \prod_{p_n>R} \left(1-\frac{Rt}{p_n}\right)^2 = R^{-\nu}\bar{w}_{X,\nu,R}(Rt), \qquad t\in [0,1].
\end{equation}

By Lemma \ref{lem:rescalingweights}, we then know that the associated normalized Christoffel-Darboux kernel transforms as
\[K_{N(R)}(x,y; \bar{w}_{X,\nu,R})=\frac{1}{R} K_{N(R)}(\tfrac{x}{R},\tfrac{y}{R}; w_{X,\nu,R}).\]
Hence, in order to prove that Theorem \ref{thm:maintheorem} holds, it suffices to prove that we have
\begin{equation}
\label{eq:rescaledmaintheorem}
\lim_{R\rightarrow \infty} \frac{1}{R} K_{N(R)}(\tfrac{x}{R},\tfrac{y}{R}; w_{X,\nu,R}) = J_\nu(x,y),
\end{equation}
uniformly on compact sets. Readers familiar with the Bessel kernel appearing at (for example) the hard edge of certain random matrix ensembles will recognize this kind of asymptotic result when properly rescaling and `zooming in on the hard edge'. For this, we note that $N(R)\rightarrow \infty$ as $R\rightarrow \infty$.

We prove \eqref{eq:rescaledmaintheorem} with the technique used by Kuijlaars and Mi\~{n}a-D\'{i}az in \cite{Kuijlaars_MinaDiaz}. Namely, we approximate the weight $w_{X,\nu,R}$ by exponential weights in Section \ref{subsec:estimatingweights}. Subsequently, we discuss the asymptotics of the kernels associated to these approximating weights in Section \ref{subsec:asymptoticsapprox}. Due to the technical nature of the associated Riemann-Hilbert problem, we postpone the proof of these asymptotics until Section \ref{sec:RHproblem}. First we will show, in Section \ref{subsec:asymptoticsactualweight},  that the asymptotics of the approximating weights actually imply \eqref{eq:rescaledmaintheorem} and hence Theorem \ref{thm:maintheorem}. For this, we use techniques introduced by Lubinsky in \cite{Lubinsky} that are also used in \cite{Kuijlaars_MinaDiaz}.

\subsection{Approximating weights}
\label{subsec:estimatingweights}
As mentioned above, we approximate the weight $w_{X,\nu,R}$ by exponential weights; such weights are well-studied in the literature.  For this, we use results from \cite{Kuijlaars_MinaDiaz}; we combine these results in the following Lemma.

\begin{lemma}[Kuijlaars{ }-{ }Mi\~{n}a-D\'{i}az, \cite{Kuijlaars_MinaDiaz}] \label{lem:approxKuijMin}
Suppose that $(q_n)_{n\in \mathbb{Z}}$ is a strictly increasing doubly infinite sequence such that the following requirements hold:
\begin{itemize}
\item[(a)] The points are indexed such that
\[\cdots < q_{-2} < q_{-1} < 0 \leq q_0< q_1 < \cdots,\]
\item[(b)] the limit
\[\lim_{S\rightarrow \infty} \sum_{0<\lvert q_n \rvert < S} \frac{1}{q_n}\]
exists,
\item[(c)] and
\[\lim_{n \rightarrow \pm \infty} \frac{q_n}{n} =1.\]
\end{itemize}
Furthermore, let $\tilde{N}(R)$ be an integer depending on $R>0$ in such a way that
\begin{equation}
\label{eq:ratioNRsine}
\lim_{R\rightarrow \infty} \frac{\tilde{N}(R)}{R}=2,
\end{equation}
and define
\begin{equation}
\label{eq:epsilonforsine}
\epsilon_R = \frac{2R}{\tilde{N}(R)} \sum_{\lvert q_n \rvert >R} \frac{1}{q_n}.
\end{equation}
Lastly, let
\begin{equation}
\label{eq:weightsine}
\tilde{w}_R(t)= \prod_{\lvert q_n\rvert >R} \Big(1-\frac{Rt}{q_n}\Big)^2, \qquad t\in [-1,1],
\end{equation}
and
\begin{equation}
\label{eq:sineexternalfield}
\tilde{V}(t)=(1+t)\log(1+t) + (1-t)\log(1-t), \qquad t\in [-1,1].
\end{equation}
Then the following two approximations hold:
\begin{enumerate}
\item For every $\alpha>1$, there is an $R_\alpha>0$ such that if $R\geq R_\alpha$, then
\begin{equation}
\label{eq:approxsineabove}
\tilde{w}_R(t) \leq \exp(\tilde{N}(R)(\tilde{V}(\tfrac{t}{\alpha})+\epsilon_R t)), \qquad t\in [-1,1].
\end{equation}
\item For every $\alpha>1$ and $\beta\in (0,1)$, there is an $R_{\alpha,\beta}>0$ such that $R\geq R_{\alpha,\beta}$ implies
\begin{equation}
\label{eq:approxsinebelow}
\tilde{w}_R(t) \geq \exp(-\tilde{N}(R)(\tilde{V}(\alpha t)+\epsilon_R t)), \qquad t\in [-\beta,\beta].
\end{equation}
\end{enumerate} 
\end{lemma}

These approximations were used to study the conditional measures of the sine process. Surprisingly, although we are studying the Bessel point process here, we can immediately use these approximations, by making a suitable transformation. Instead of the external field \eqref{eq:sineexternalfield}, we consider 
\begin{equation}
\label{eq:externalfield}
V(t) = 2(1+\sqrt{t}) \log(1+\sqrt{t}) + 2(1-\sqrt{t}) \log(1-\sqrt{t}), \qquad t\in [0,1].
\end{equation} 
Clearly, we have 
\begin{equation}
\label{eq:relationpotentialsBesselSine}
V(t)=2\tilde{V}(\sqrt{t}),
\end{equation}
where $\tilde{V}$ is given by \eqref{eq:sineexternalfield}. Then, we find the following approximations.

\begin{proposition}
\label{prop:approx}
Suppose that $X=(p_n)_{n=1}^\infty$ is an increasing sequence of positive numbers satisfying the growth condition \eqref{eq:growthconditionpoints}. Then, for every $\gamma>1$ there exists an $R_\gamma>0$ such that $R\geq R_\gamma$ implies
\begin{equation}
\label{eq:approx}
t^{\nu}  \exp(-N(R) V(\gamma t)) \mathfrak{1}_{[0,\gamma^{-2}]}(t) \leq w_{X,\nu,R}(t)
\leq t^{\nu}  \exp(-N(R) V(t/\gamma)), \qquad t\in [0,1],
\end{equation}
where $V$ is defined by \eqref{eq:externalfield} and $N(R)$ is the number of points in $X\cap (0,R]$.
\end{proposition}
\begin{proof}
We notice that
\begin{equation}
\label{eq:splitweightintwoparts}
t^{-\nu} w_{X,\nu,R}(t) = \prod_{p_n > R} \left(1-\frac{R t}{p_n}\right)^2
= \prod_{p_n>R} \left(1-\frac{\frac{\sqrt{R}}{\pi} \sqrt{t}}{\frac{\sqrt{p_n}}{\pi}}\right)^2 \left(1+\frac{\frac{\sqrt{R}}{\pi} \sqrt{t}}{\frac{\sqrt{p_n}}{\pi}}\right)^2, \qquad t\in [0,1].
\end{equation}
Motivated by this, we define the doubly infinite sequence $(q_n)_{n\in \mathbb{Z}}$ by
\[q_n = \begin{cases}
\frac{\sqrt{p_n}}{\pi} &\text{if $n\geq 1$}\\
-\frac{\sqrt{p_{1-n}}}{\pi} &\text{if $n\leq 0$}.
\end{cases}\]
We note that this sequence $(q_n)_{n\in \mathbb{Z}}$ satisfies the requirements (a), (b) of Lemma \ref{lem:approxKuijMin} by construction, and (c) due to the growth condition \eqref{eq:growthconditionpoints} on $X$. We define $\tilde{N}(R)=2N(\pi^2R^2)$. Since we have that
\begin{equation}
\label{eq:RandNR}
\lim_{R\rightarrow \infty} \frac{R}{N(R)^2}=\pi^2
\end{equation}
by \eqref{eq:growthconditionpoints}, we also have that \eqref{eq:ratioNRsine} holds. Furthermore, we have that  $q_n=-q_{1-n}$ for all $n\in \mathbb{Z}$, so \eqref{eq:epsilonforsine} gives us that $\epsilon_R=0$ for all $R>0$. The factorization in \eqref{eq:splitweightintwoparts} is especially useful since it implies that
\begin{equation}
\label{eq:wandwtilde}
t^{-\nu} w_{X,\nu,R}(t) = \prod_{\lvert q_n \rvert > \frac{\sqrt{R}}{\pi}} \left(1-\frac{\frac{\sqrt{R}}{\pi} \sqrt{t}}{q_n}\right)^2 = \tilde{w}_{\frac{\sqrt{R}}{\pi}}(\sqrt{t}), \qquad t\in [0,1],
\end{equation}
where $\tilde{w}$ is given by \eqref{eq:weightsine}. Indeed, this allows us to transfer the approximations in \eqref{eq:approxsineabove}- \eqref{eq:approxsinebelow} to the approximations of the weight $w_{X,\nu,R}$.

Namely, let $\gamma>1$ and let $\tilde{w}_R$ and $\tilde{V}$ be defined as in \eqref{eq:weightsine} and \eqref{eq:sineexternalfield}, respectively. Then since $\sqrt{\gamma}>1$, there is an $R_1>0$ such that $R\geq R_1$ implies that
\begin{equation*}
\tilde{w}_R(t) \leq \exp\left(\tilde{N}(R)\tilde{V}\left(\tfrac{t}{\sqrt{\gamma}}\right)\right), \qquad t\in [-1,1],
\end{equation*}
by using \eqref{eq:approxsineabove}. Furthermore, if we take $\alpha=\sqrt{\gamma}$ and $\beta=\gamma^{-1}$, we have that $\alpha>1$ and $\beta\in (0,1)$, so if we apply \eqref{eq:approxsinebelow} we obtain that there is an $R_2>0$ such that $R\geq R_2$ implies that 
\begin{equation*}
\tilde{w}_R(t) \geq \mathfrak{1}_{[0,\gamma^{-1}]}(t) \exp(-\tilde{N}(R)\tilde{V}(\sqrt{\gamma} t)), \qquad t\in [-1,1].
\end{equation*}

Now define $R_\gamma = \pi^2 \max(R_1,R_2)^2$. Then $R\geq R_\gamma$ implies that $\frac{\sqrt{R}}{\pi}\geq \max(R_1,R_2)$, and that implies that we have
\[\mathfrak{1}_{[0,\gamma^{-1}]}(\sqrt{t}) \exp(-\tilde{N}(\sqrt{R}/\pi)\tilde{V}(\sqrt{\gamma t})) \leq \tilde{w}_{\frac{\sqrt{R}}{c}}(\sqrt{t}) \leq \exp(\tilde{N}(\sqrt{R}/\pi) \tilde{V}(\sqrt{t/\gamma})), \qquad t\in [0,1].\]
Now using the fact that $\tilde{N}(\sqrt{R}/\pi)=2N(R)$ and $V(t)=2\tilde{V}(\sqrt{t})$, combined with \eqref{eq:wandwtilde}, gives us
\[t^\nu \mathfrak{1}_{[0,\gamma^{-2}]}(t) \exp(-N(R) V(\gamma t)) \leq w_{X,\nu,R}(t) \leq t^\nu \exp(N(R) V(t/\gamma)), \qquad t\in [0,1],\]
for every $R\geq R_\gamma$. The existence of such an $R_\gamma$ for every $\gamma>1$ is exactly what we set out to prove.
\end{proof}

We refer to the weights which give a lower and an upper bound for the weight $w_{X,\nu,R}$ in \eqref{eq:approx} as the \textbf{approximating weights}. We reserve a notation for these weights, namely
\begin{align}
\label{eq:varyingexponentialweightplus}
\omega_{\gamma,n,\nu}^+ (t)&= t^\nu \exp(-nV(t/\gamma)), &\qquad t\in [0,1], \\
\label{eq:varyingexponentialweightminus} 
\omega_{\gamma,n,\nu}^- (t)&= t^\nu \mathfrak{1}_{[0,\gamma^{-2}]}(t) \exp(-n V(\gamma t)), &\qquad t\in [0,1],
\end{align}
where, as before, $V$ is defined by \eqref{eq:externalfield}, $\gamma>1$, $n\geq 1$ is an integer and $\nu>-1$. Using this notation, we have that \eqref{eq:approx} precisely becomes
\begin{equation}
\label{eq:approxinnotation}
\omega_{\gamma,N(R),\nu}^-(t) \leq w_{X,\nu,R}(t) \leq \omega_{\gamma,N(R),\nu}^+(t).
\end{equation} 
We note that the two weights in \eqref{eq:varyingexponentialweightplus}-\eqref{eq:varyingexponentialweightminus} are transformed to each other according to
\begin{equation}
\label{eq:transformplusminus}
\omega_{\gamma,n,\nu}^-(t)=\frac{1}{\gamma^{2\nu}} \omega_{\gamma,n,\nu}^+(\gamma^2 t), \qquad t\in [0,\gamma^{-2}].
\end{equation}

\subsection{The relevant equilibrium measure}
\label{subsec:equilibrium}
As mentioned above, we do not study the asymptotics of the weight $w_{X,\nu,R}$ directly, but instead study the asymptotics of the approximating weights $\omega_{\gamma,N(R),\nu}^{\pm}$. In fact, by \eqref{eq:transformplusminus}, we are only required to study the weight $\omega_{\gamma,n,\nu}^+$. Since this weight is of the form \eqref{eq:varyingexponentialweightplus}, a reader familiar with this kind of asymptotics knows that one is interested in the \textit{equilibrium measure} of the external field
\begin{equation}
\label{eq:Vgamma}
V_\gamma(t) =V(t/\gamma), \qquad t\in [0,1].
\end{equation}
Therefore, we compute this equilibrium measure before turning to the actual asymptotics.

For an arbitrary external field $\tilde{V}$ on some interval $[a,b]$, the associated \textbf{equilibrium measure} $\mu_{\tilde{V}}$ is the unique probability measure $\mu_{\tilde{V}}$ with support in $[a,b]$ for which there is a constant $\ell_{\tilde{V}}$ such that the following equation holds:
\begin{equation}
\label{eq:equilibriumproblem}
2\int_0^1 \log|x-s| d\mu_{\tilde{V}}(s)  
\left\{\begin{array}{ll} 
= \tilde{V}(x) + \ell_{\tilde{V}},& x\in \operatorname{supp} \mu_{\tilde V},\\
\leq \tilde{V}(x) + \ell_{\tilde{V}},& x\in [a,b]\setminus \operatorname{supp} \mu_{\tilde V}.
\end{array}\right.
\end{equation}
The search for this equilibrium measure is often referred to as the \textit{equilibrium problem}. We first solve this equilibrium problem on $[0,1]$ for the external field $V$, that was defined in \eqref{eq:externalfield}.
\begin{lemma}
\label{lem:eqmeasureforexternalfield}
The equilibrium measure on $[0,1]$ for the external field $V$ given by \eqref{eq:externalfield} is 
\begin{equation}
\label{eq:equilibriumgamma1}
\frac{d\mu_V(s)}{ds}=\frac{1}{2\sqrt{s}}, \qquad s\in [0,1],
\end{equation}
and the constant $\ell_V=0$.
\end{lemma}
\begin{proof}
For any $s\in [0,1]$ we have by making the transformation $s=t^2$
\begin{align*}
2\int_0^1 \log|x-s| \frac{1}{2\sqrt{s}} ds &= 2 \int_0^1 \log \lvert x-t^2 \rvert dt \\
&=2 \int_0^1 (\log \lvert \sqrt{x} +t \rvert +\log \lvert \sqrt{x} -t \rvert) dt \\
&=2 \int_{-1}^1 \log \lvert \sqrt{x} -t \rvert dt \\
&= V(x),
\end{align*}
where in the last step we used the computation of the equilibrium measure of the external field that plays a role in \cite{Kuijlaars_MinaDiaz}.
\end{proof}

Now we turn to the external field $V_\gamma$, as in \eqref{eq:Vgamma}. For brevity of notation, we denote the equilibrium measure on $[0,1]$ associated to this deformed external field by $\mu_\gamma$. We have an explicit expression.
\begin{lemma}
Let $\gamma>1$. The equilibrium measure $\mu_\gamma$ on $[0,1]$ for the external field $V_\gamma$ has a density with support $[0,1]$, given explicitly by
\begin{align}
\label{eq:equilibriummeasuregamma}
\frac{d\mu_\gamma}{ds} &= \frac{1}{\sqrt{\gamma s}} \left(\frac{1}{2}+\frac{1}{\pi} \sqrt{\frac{\gamma-1}{1-s}}
-\frac{1}{\pi} \arctan\left(\sqrt{\frac{\gamma-1}{1-s}}\right)\right), & s\in (0,1).
\end{align}
\end{lemma}
\begin{proof}
First we note that the measure $m_\gamma$ with density
\begin{equation}
\frac{dm_\gamma}{ds}(s) := \frac{1}{2\sqrt{\gamma s}}, \qquad s\in[0,\gamma],
\end{equation}
is a probability measure on $[0,\gamma]$, which can easily be checked. Next, we note that for any $x\in[0,\gamma]$ we have that
\begin{equation*}
2\int_0^\gamma \log|x-s| dm_\gamma(s) = \int_0^1 \log\lvert x-\gamma t \rvert \frac{dt}{2\sqrt{t}} = 2\log(\gamma) + 2\int_0^1 \log \Big\lvert \tfrac{x}{\gamma} -t \Big\rvert \frac{dt}{2\sqrt{t}}.
\end{equation*} 
Now, making use of Lemma \ref{lem:eqmeasureforexternalfield}, we see that the last term is in fact equal to $V_\gamma(x)$. Hence we obtain that
\begin{equation*}
2\int_0^\gamma \log|x-s| dm_\gamma(s) = V_\gamma(x)+2\log(\gamma),
\end{equation*}
which means that the measure $m_\gamma$ satisfies all the requirements of the equilibrium problem \eqref{eq:equilibriumproblem}, except that it is supported on $[0,\gamma]$ instead of a subset of $[0,1]$. There is a general framework for dealing with this peculiarity, namely via the use of \textit{balayages}. According to equation (4.47) in \cite[Chapter II]{Saff_Totik}, we now have for $s\in [0,1]$
\begin{align*}
\frac{d\mu_\gamma}{ds} &= \frac{dm_\gamma}{ds}+\frac{1}{\pi} \int_1^\gamma \frac{\sqrt{x(x-1)}}{(x-s)\sqrt{s(1-s)}} dm_\gamma(x)\\
&= \frac{1}{2\sqrt{\gamma s}} + \frac{1}{2\pi\sqrt \gamma} \frac{1}{\sqrt{s(1-s)}}\int_1^\gamma \frac{\sqrt{x-1}}{x-s} dx.
\end{align*}
It can be checked that for $s\in(0,1)$
\[
\int_1^\gamma \frac{\sqrt{x-1}}{x-s} dx
= 2\sqrt{\gamma-1} - 2\sqrt{1-s} \arctan\left(\sqrt{\frac{\gamma-1}{1-s}}\right),
\]
and therefore \eqref{eq:equilibriummeasuregamma} follows.
\end{proof}
We notice that the behaviour of the equilibrium measure around $s=0$ and $s=1$ immediately follows from \eqref{eq:equilibriummeasuregamma}. Namely, for fixed $\gamma>1$ we have
\begin{align}
\frac{d\mu_\gamma}{ds} &= s^{-\frac{1}{2}} \frac{1}{\sqrt \gamma} \left(\frac{1}{2} + \frac{\sqrt{\gamma-1}-\arctan{\sqrt{\gamma-1}}}{\pi}\right)+\mathcal O(s^\frac{1}{2}), & s\to 0^+ \label{eq:equilbriumaround0}\\
\frac{d\mu_\gamma}{ds} &= \frac{1}{\pi} \sqrt{\frac{\gamma-1}{\gamma}} (1-s)^{-\frac{1}{2}}+\mathcal O((1-s)^\frac{1}{2}), & s\to 1^-.\label{eq:equilbriumaround1}
\end{align}

\subsection{Asymptotics for the approximating weights}
\label{subsec:asymptoticsapprox}
For the asymptotics of the approximating weight $\omega_{\gamma,n,\nu}^+$, we use the constant
\begin{equation}
\label{eq:defcgamma}
c_\gamma=\frac{\pi^2 \gamma}{(\pi+2(\sqrt{\gamma-1}-\arctan \sqrt{\gamma-1}))^2},
\end{equation} 
that we consider for $\gamma>1$. The reason for this specific choice of constant follows from the Riemann-Hilbert problem that we study in  Section \ref{sec:RHproblem}. We note that $c_\gamma>1$ and that
\begin{equation}
\lim_{\gamma \rightarrow 1^+} c_\gamma =1.
\end{equation} 
We then have the following.
\begin{proposition}
\label{prop:resultRHprob}
Suppose that $\gamma>1$ and let $c_\gamma$ be defined by \eqref{eq:defcgamma}. Then the normalized Christoffel-Darboux kernel of the weight in \eqref{eq:varyingexponentialweightplus} satisfies
\begin{equation}
\label{eq:resultRHprob}
\lim_{n\rightarrow \infty} \frac{c_\gamma}{\pi^2 n^2} K_n\left(\frac{c_\gamma x}{\pi^2 n^2},\frac{c_\gamma y}{\pi^2 n^2}; \omega_{\gamma,n,\nu}^+\right) =  J_\nu(x,y),
\end{equation}
uniformly for $(x,y)$ in compact subsets of $(0,\infty)^2$.
\end{proposition}
Proposition \ref{prop:resultRHprob} can be proven using standard Riemann-Hilbert techniques. These techniques were developed as a tool to study the asymptotics of orthogonal polynomials. This started with the seminal paper by Fokas, Its and Kitaev \cite{Fokas_Its_Kitaev} who studied orthogonal polynomials with respect to weights on contours in the complex plane. These techniques were refined by (amongst others) Deift, Kriecherbauer, McLaughlin, Venakides and Zhou \cite{Deift,Deift_Kriecherbauer_McLaughlin_Venakides_Zhou} for weights on the real line.

The use of Riemann-Hilbert problems was extended to weights that were not supported on the real line, but on bounded intervals. Kuijlaars, McLaughlin, Van Assche and Vanlessen \cite{Kuijlaars_2004,Kuijlaars_2002} solved these problems in a way that has since become the standard. The main difference between the analysis on a bounded interval instead of the whole real line lies in the analysis that needs to be undertaken at the end-points. 

In a Riemann-Hilbert problem, the local parametrices are meant to deal with this local analysis. In our case, we are interested in the local parametrix around $z=0$, where the equilibrium measure blows up as an inverse square root, according to \eqref{eq:equilbriumaround0}. It is well-known that if this is the case and the weight is of the form \eqref{eq:varyingexponentialweightplus}, i.e., of varying exponential type, then one can use a modification of the \textit{Bessel parametrix} that was defined in \cite{Kuijlaars_2004} to find a solution of the local parametrix problem at hand. This is precisely what we do in Section \ref{sec:RHproblem}. Similar analyses have been carried out in for example \cite[Section 7.1]{Celsus_Silva}, \cite[Section IV.D]{Charlier_Claeys}, and \cite[Section 4.1.4]{Deano}. To explain the appearance of the constant $c_\gamma$ in \eqref{eq:resultRHprob}, we provide the details of the Riemann-Hilbert problem in Section \ref{sec:RHproblem}. The reader who is familiar with Riemann-Hilbert problems will already be familiar with these details.

\begin{remark}
If one would take $\gamma=1$, one obtains the equilibrium measure in \eqref{eq:equilibriumgamma1}. Note specifically that its behaviour around $s=1$ is qualitatively different from \eqref{eq:equilbriumaround1}; it does not blow up as $(1-s)^{-\frac{1}{2}}$ but in fact has leading order $O(1)$. We would not know how to solve the corresponding local parametrix problem.
\end{remark}

By the transformative property \eqref{eq:transformplusminus} of the approximating weights and the general transformation rule \eqref{eq:norCDkernrescaled} for normalized Christoffel-Darboux kernels, Proposition \ref{prop:resultRHprob} immediately implies the asymptotics for the other approximating weight.
\begin{corollary}
\label{cor:resultRHprobminus}
The normalized Christoffel-Darboux kernel of the weight in \eqref{eq:varyingexponentialweightminus} satisfies
\begin{equation}
\label{eq:resultRHprobminus}
\lim_{n\rightarrow \infty} \frac{1}{c_\gamma \pi^2 n^2} K_n\left(\frac{x}{c_\gamma \pi^2 n^2},\frac{y}{c_\gamma \pi^2 n^2}; \omega_{\gamma,n,\nu}^-\right) = J_\nu(x,y),
\end{equation}
uniformly for $(x,y)$ in compact subsets of $(0,\infty)^2$, for every $\gamma>1$.
\end{corollary}
For further analysis, it is also important to note that the approximating weights satisfy the following asymptotic behaviour uniformly on compact sets:
\begin{equation}
\label{eq:asymptoticapproxweight}
\lim_{n\rightarrow \infty} n^{2\nu} \omega_{\gamma,n,\nu}^\pm\left(\frac{x}{n^2}\right) = x^\nu.
\end{equation}
This follows readily from the behavior $V(x)=\mathcal O(x)$ as $x\to 0$. We especially note that the limit is independent of $\gamma$.
We proceed by considering the non-normalized kernels of the approximating weights.

\begin{corollary}
\label{cor:nonnormalizedresult}
The non-normalized Christoffel-Darboux kernels of the weights in \eqref{eq:varyingexponentialweightplus} and \eqref{eq:varyingexponentialweightminus} have the following asymptotic behaviour:
\begin{align}
\lim_{n \rightarrow \infty} \frac{1}{(\pi n)^{2+2\nu}} \widehat{K}_n\left(\frac{x}{\pi^2 n^2},\frac{y}{\pi^2 n^2}; \omega_{\gamma,n,\nu}^+\right) &= (xy)^{-\nu/2} \frac{1}{c_\gamma} J_\nu\left(\frac{x}{c_\gamma},\frac{y}{c_\gamma}\right), \label{eq:resultRHprobhatplus} \\
\lim_{n \rightarrow \infty} \frac{1}{(\pi n)^{2+2\nu}} \widehat{K}_n\left(\frac{x}{\pi^2 n^2},\frac{y}{\pi^2 n^2}; \omega_{\gamma,n,\nu}^-\right) &= (xy)^{-\nu/2} c_\gamma J_\nu\left(c_\gamma x,c_\gamma y\right), \label{eq:resultRHprobhatminus}
\end{align}
uniformly for $(x,y)$ in compact subsets of $(0,\infty)^2$, for all $\gamma>1$.
\end{corollary}
\begin{proof}
This follows directly by replacing $(x,y)$ by $c_\gamma^{\mp 1} (x,y)$ and combining \eqref{eq:resultRHprob} and \eqref{eq:resultRHprobminus} with the relationship between the normalized and non-normalized Christoffel-Darboux kernel and the asymptotic behaviour \eqref{eq:asymptoticapproxweight}.
\end{proof}

\subsection{Asymptotics for the actual weight}
\label{subsec:asymptoticsactualweight}
In this section, we use the above asymptotics for the kernels of the approximating weights $\omega_{\gamma,n,\nu}^\pm$ to obtain the asymptotics for the kernel of the weight $w_{X,\nu,R}$, i.e. Theorem \ref{thm:maintheorem}. Informally, we let $\gamma \rightarrow 1^+$, and by that `squeeze in' the kernel of interest by using \eqref{eq:approxinnotation}. 

\begin{proposition}
\label{prop:maintheoremforNRhat}
Suppose that $X=(p_n)_{n=1}^\infty$ is a strictly increasing sequence of positive numbers satisfying the growth condition \eqref{eq:growthconditionpoints}, and for every $R>0$, let $N(R)$ be the number of points in $X\cap (0,R]$. Then, uniformly for $(x,y)$ in compact subsets of $(0,\infty)^2$, we have that
\begin{equation}
\label{eq:maintheoremforNRhat}
\lim_{R \rightarrow \infty} \frac{1}{(\pi N(R))^{2+2\nu}} \widehat{K}_{N(R)}\left(\frac{x}{\pi^2 N(R)^2},\frac{y}{\pi^2 N(R)^2}; w_{X,\nu,R}\right) = (xy)^{-\frac{\nu}{2}} J_\nu(x,y).
\end{equation}
\end{proposition}
\begin{proof}
We assume Proposition \ref{prop:resultRHprob} and in particular Corollary \ref{cor:nonnormalizedresult}. First we prove \eqref{eq:maintheoremforNRhat} on the diagonal, that is $x=y$, and then we extend this result to all $(x,y)$ by using a technique by Lubinsky, developed in \cite{Lubinsky}.

Now let $\gamma>1$. By Proposition \ref{prop:approx}, we know that there is an $R_\gamma>0$ such that for all $R\geq R_\gamma$, we have that
\begin{align} \label{eq:weightsIneqLub}
\omega_{\gamma,N(R),\nu}^-(t) \leq w_{X,\nu,R}(t) \leq \omega_{\gamma,N(R),\nu}^+(t), \qquad t\in [0,1].
\end{align}
Then, by using the extremal property of the Christoffel function (cf. \eqref{eq:minpropChr}), we also have for all $R\geq R_\gamma$ and for all $x\in [0,1]$, that
\begin{align} \label{eq:kernelsIneqLub}
\widehat{K}_{N(R)}(x,x; \omega_{\gamma,N(R),\nu}^+) \leq \widehat{K}_{N(R)}(x,x; w_{X,\nu,R}) \leq \widehat{K}_{N(R)}(x,x; \omega_{\gamma,N(R),\nu}^-).
\end{align}
Due to the differentiability of $(x,y)\mapsto \sqrt{x y} J_\nu(x,y)$ on $(0,\infty)^2$ we have uniformly for $(x,y)$ in any compact subset of $(0,\infty)^2$ that there exists a constant $M>0$ such that
\begin{equation}
\label{eq:existenceM1}
(xy)^{-\nu/2} \Big\lvert J_\nu(x,y) - \frac{1}{c_\gamma} J_\nu\left(\frac{x}{c_\gamma},\frac{y}{c_\gamma}\right) \Big\rvert \leq M\left(1-\frac{1}{c_\gamma}\right).
\end{equation}
Now let $x$ be in a compact set $S$ and take $R$ big enough such that $x_R:=\frac{x}{\pi^2 N(R)^2}\in (0,1)$ for all $x\in S$. Furthermore, let $\varepsilon>0$. By \eqref{eq:kernelsIneqLub} and \eqref{eq:existenceM1} we can find a $\gamma>1$ such that uniformly on $S$
\begin{align} \label{eq:+K-1}
\frac{1}{(\pi N(R))^{2+2\nu}} &\widehat{K}_{N(R)}(x_R,x_R; \omega_{\gamma,N(R),\nu}^+) - x^{-\nu} \frac{1}{c_\gamma} J_\nu\left(\frac{x}{c_\gamma},\frac{x}{c_\gamma}\right)  - \varepsilon \\ \label{eq:+K-2}
&\leq \frac{1}{(\pi N(R))^{2+2\nu}} \widehat{K}_{N(R)}(x_R,x_R; w_{X,\nu,R}) - x^{-\nu} J_\nu(x,x) \\ \label{eq:+K-3}
&\leq \frac{1}{(\pi N(R))^{2+2\nu}} \widehat{K}_{N(R)}(x_R,x_R; \omega_{\gamma,N(R),\nu}^-) - x^{-\nu} c_\gamma J_\nu\left(c_\gamma x,c_\gamma x\right) + \varepsilon.
\end{align}
Here we have tacitly used that $c_\gamma$ tends to $1$ as $\gamma\to 1^+$. Then using the uniform convergence in Corollary \ref{cor:nonnormalizedresult} for the lower bound \eqref{eq:+K-1} and the upper bound \eqref{eq:+K-3} we infer that for $R$ big enough we have that uniformly for $x\in S$
\begin{align} \label{eq:eKe}
 -2\varepsilon \leq \frac{1}{(\pi N(R))^{2+2\nu}} \widehat{K}_{N(R)}(x_R,x_R; w_{X,\nu,R}) - x^{-\nu} J_\nu(x,x) \leq 2\varepsilon.
\end{align}
This proves \eqref{eq:maintheoremforNRhat} on the diagonal. 

For general $x,y>0$ we can use the result on the diagonal to prove that \eqref{eq:maintheoremforNRhat} holds on compact subsets. By \eqref{eq:weightsIneqLub} we may use an inequality by Lubinsky \cite{Lubinsky}, given by
\begin{multline} \label{eq:LubOmega+w}
\Big \lvert \widehat{K}_{N(R)}(x_R,y_R; w_{X,\nu,R}) - \widehat{K}_{N(R)}(x_R, y_R; \omega_{\gamma,N(R),\nu}^+) \Big \rvert^2 \\
 \leq \Big\lvert \widehat{K}_{N(R)}(y_R,y_R; w_{X,\nu,R}) \Big\rvert \left(\widehat{K}_{N(R)}(x_R,x_R; w_{X,\nu,R})-\widehat{K}_{N(R)}(x_R,x_R; \omega_{\gamma,N(R),\nu}^+)\right).
\end{multline}
The proposition now follows by a similar argument that we applied for the diagonal, using the uniform convergence of \eqref{eq:maintheoremforNRhat} on compact subsets of the diagonal that we just proved, \eqref{eq:resultRHprobhatplus} and \eqref{eq:existenceM1}. 
\end{proof}

By combining \eqref{eq:approxinnotation} and \eqref{eq:asymptoticapproxweight}, we immediately see that we have
\begin{equation}
\label{eq:asympforweight}
\lim_{R\rightarrow \infty} \pi^{2\nu} N(R)^{2\nu} w_{X,\nu,R}\left(\frac{x}{\pi^2 N(R)^2}\right) = x^\nu,
\end{equation} 
uniformly on compact subsets. From this, we can directly conclude the limiting behaviour of the normalized Christoffel-Darboux kernel of the weight $w_{X,\nu,R}$.

\begin{corollary}
\label{cor:maintheoremforNR}
Suppose that $X=(p_n)_{n=1}^\infty$ is a strictly increasing sequence of positive numbers satisfying the growth condition \eqref{eq:growthconditionpoints}, and for every $R>0$, let $N(R)$ be the number of points in $X\cap [0,R]$. Then, uniformly for $(x,y)$ in compact subsets of $(0,\infty)^2$, we have that
\begin{equation}
\label{eq:maintheoremforNR}
\lim_{R \rightarrow \infty} \frac{1}{\pi^2 N(R)^2} K_{N(R)}\left(\frac{x}{\pi^2 N(R)^2},\frac{y}{\pi^2 N(R)^2}; w_{X,\nu,R}\right) =  J_\nu(x,y).
\end{equation}
\end{corollary}
\begin{proof}
Combine the relationship between the normalized and the non-normalized Christoffel-Darboux kernel with Proposition \ref{prop:maintheoremforNRhat} and \eqref{eq:asympforweight}.
\end{proof}

As a corollary, we have our main theorem.

\begin{proof}[Proof of Theorem \ref{thm:maintheorem}]
This directly follows from Corollary \ref{cor:maintheoremforNR} and \eqref{eq:RandNR}.
\end{proof}

%% file: shortRHproblem.tex
\section{The Riemann-Hilbert problem}
\label{sec:RHproblem}
We now give the details for the RH-problem for the weight $\omega_{\gamma,n,\nu}^+$, which proves Proposition \ref{prop:resultRHprob}. For this, we fix $\gamma>1$ and $\nu>-1$. We also already note that $\omega_{\gamma,n,\nu}^+(x)=O(x^\nu)$ as $x\downarrow 0$ and that $\omega_{\gamma,n,\nu}^+$ is bounded as $x\to 1$, uniformly in $n$. Furthermore, for convenience, we define the function
\begin{equation}
\label{eq:defhnu}
h_\nu(z) = \left\{\begin{array}{rl} 1 & \nu>0\\ \log z & \nu=0\\ z^\nu & \nu <0. \end{array}\right.
\end{equation}

The Riemann-Hilbert problem that we consider is the following:
\begin{itemize}
\item[RH-Y1] $Y : \mathbb{C}\setminus [0,1]\to \mathbb{C}^{2\times 2}$ is analytic.
\item[RH-Y2] For $x\in(0,1)$ we have (oriented away from the origin)
\[Y_+(x) = Y_-(x) \begin{pmatrix} 1 & \omega_{\gamma,n,\nu}^+(x)\\ 0 & 1\end{pmatrix}.\]
\item[RH-Y3] As $z\to \infty$ we have
\[Y(z) = (I+O(1/z)) \begin{pmatrix} z^n & 0\\ 0 & z^{-n}\end{pmatrix}.\]
\item[RH-Y4a] As $z\to 0$ we have
\[Y(z) = O \begin{pmatrix} 1 & h_\nu(z)\\ 1 & h_\nu(z)\end{pmatrix}.\]
\item[RH-Y4b] As $z\to 1$ we have
\[Y(z) = O \begin{pmatrix} 1 & \log(z-1)\\ 1 & \log(z-1) \end{pmatrix}.\]
\end{itemize}
In \cite{Kuijlaars_2004}, this Riemann-Hilbert problem was used for the first time to study the asymptotics of orthogonal polynomials on a bounded interval. We explicitly state that its solution is given by
\begin{equation}
\label{eq:explicitY}
Y(z) = 
\begin{pmatrix}
\gamma_n^{-1} \varphi_n(z;\omega_{\gamma,n,\nu}^+) & \frac{\gamma_n^{-1}}{2\pi i}  \int_0^1 \frac{\varphi_n(s;\omega_{\gamma,n,\nu}^+)\omega_{\gamma,n,\nu}^+(s)}{s-z} ds\\
-2\pi i\gamma_{n-1} \varphi_{n-1}(z;\omega_{\gamma,n,\nu}^+) & -\gamma_{n-1}  \int_0^1 \frac{\varphi_{n-1}(s;\omega_{\gamma,n,\nu}^+) \omega_{\gamma,n,\nu}^+(s)}{s-z} ds
\end{pmatrix},
\end{equation}
where the constant $\gamma_n$ is the leading coefficient of the orthonormal polynomial $\varphi_n(z;\omega_{\gamma,n,\nu}^+)$. The normalized Christoffel-Darboux kernel is explicitly given in terms of the solution $Y$ by:
\begin{equation}
\label{eq:kernelinY}
K_n(x,y;\omega_{\gamma,n,\nu}^+) = \frac{1}{2\pi i (x-y)} \sqrt{\omega_{\gamma,n,\nu}^+(x) \omega_{\gamma,n,\nu}^+(y)} 
\begin{pmatrix} 0 & 1\end{pmatrix}
Y_+(y)^{-1} Y_+(x) 
\begin{pmatrix}
1\\ 0
\end{pmatrix}.
\end{equation}

We note that in our analysis we use the convention to take the principle branch for logarithms and power functions, i.e. these will have cut $(-\infty,0]$ and they will be positive for large positive values.

\subsection{First transformation: normalization}
For our fixed $\gamma>1$, we define the function $g_\gamma:\mathbb{C}\setminus (-\infty,1]\to \mathbb{C}$ by 
\begin{equation}
g_{\gamma}(z) = \int_0^1 \log(z-s) d\mu_\gamma(s),
\end{equation}
where the logarithm is defined with argument in $(-\pi,\pi)$, as usual, and $\mu_\gamma$ is the equilibrium measure \eqref{eq:equilibriummeasuregamma}. Then we have that $g_\gamma$ is analytic on $\mathbb{C}\setminus (-\infty,1]$ and it is obvious that we have 
\begin{equation} 
\label{eq:behaviourgatinfinity}
g_\gamma(z) = \log(z)+\mathcal O\left(\frac{1}{z}\right), \qquad \textrm{as} \ z\to\infty.
\end{equation} 
We note that $g_\gamma$ is bounded around $z=0$ and $z=1$, and furthermore, we have
\begin{align}
g_{\gamma,+}(x)+g_{\gamma,-}(x) &= 2 \int_0^1 \log|x-s| d\mu_\gamma(s) = V_\gamma(x)+\ell_\gamma, & x\in (0,1),\label{eq:gplus+minusinterval}\\
g_{\gamma,+}(x)-g_{\gamma,-}(x) &= 2\pi i\int_x^1 d\mu_\gamma(s), & x\in (0,1), \label{eq:gplus-minusinterval}\\
g_{\gamma,+}(x)-g_{\gamma,-}(x) &= 2\pi i, & x\leq 0 \label{eq:gplus-minusnegative}.
\end{align}
Now define the function $T: \mathbb{C}\setminus [0,1] \rightarrow \mathbb{C}^{2\times 2}$ by
\begin{equation}
\label{eq:fromYtoT}
T(z) = \begin{pmatrix} e^{-\frac{n\ell_\gamma}{2}} & 0\\ 0 & e^{\frac{n\ell_\gamma}{2}}\end{pmatrix} 
Y(z) 
\begin{pmatrix} e^{-n g_\gamma(z)} & 0\\ 0 & e^{n g_\gamma(z)}\end{pmatrix}
\begin{pmatrix} e^{\frac{n\ell_\gamma}{2}} & 0\\ 0 & e^{-\frac{n\ell_\gamma}{2}}\end{pmatrix},
\end{equation}
where $\ell_\gamma$ is the constant of the equilibrium problem \eqref{eq:equilibriumproblem} associated to $V_\gamma$. Then, by the above properties of $g$, $T$ satisfies the following Riemann-Hilbert problem.
\begin{itemize}
\item[RH-T1] $T:\mathbb{C}\setminus [0,1]\to \mathbb{C}^{2\times 2}$ is analytic.
\item[RH-T2] For $x\in(0,1)$ we have (oriented away from the origin)
\begin{equation}
\label{eq:RHT2}
T_+(x) = T_-(x) \begin{pmatrix} e^{2\pi i n \int_0^x d\mu_\gamma(s)} & x^\nu\\ 0 & e^{-2\pi i n\int_0^x d\mu_\gamma(s)}\end{pmatrix}.
\end{equation}
\item[RH-T3] As $z\to\infty$ we have $T(z) = I+O(1/z)$.
\item[RH-T4a] As $z\to 0$ we have
\[T(z) = O \begin{pmatrix} 1 & h_\nu(z)\\ 1 & h_\nu(z)\end{pmatrix}.\]
\item[RH-T4b] As $z\to 1$ we have
\[T(z) = O \begin{pmatrix} 1 & \log(z-1)\\ 1 & \log(z-1) \end{pmatrix}.\]
\end{itemize}

\subsection{Opening of the lens}
\label{subsec:openinglens}
In order to open the lens, we define the following function:
\begin{equation}
\label{eq:defphi}
\varphi_\gamma(z) =  \left\{
\begin{array}{rl} 
\displaystyle \log\left(\frac{\sqrt\gamma\sqrt{1-z}+i\sqrt z\sqrt{\gamma-1}}{\sqrt\gamma\sqrt{1-z}-i\sqrt z\sqrt{\gamma-1}}\right)
+ \sqrt{\frac{z}{\gamma}} \log\left(\frac{\sqrt{\gamma-1}+i\sqrt{1-z}}{\sqrt{\gamma-1}-i\sqrt{1-z}}\right), & \operatorname{Im}(z)>0,\\
\displaystyle \log\left(\frac{\sqrt\gamma\sqrt{1-z}-i\sqrt z\sqrt{\gamma-1}}{\sqrt\gamma\sqrt{1-z}+i\sqrt z\sqrt{\gamma-1}}\right)
+ \sqrt{\frac{z}{\gamma}} \log\left(\frac{\sqrt{\gamma-1}-i\sqrt{1-z}}{\sqrt{\gamma-1}+i\sqrt{1-z}}\right), & \operatorname{Im}(z)<0.
\end{array}
\right. 
\end{equation}
To avoid confusion, we repeat that we use the convention of taking the principal branch for the logarithm and power functions (square root functions in this case). We claim that this function $\varphi_\gamma$ is well-defined and analytic on $\mathbb{C}\setminus \mathbb{R}$. To see this, note that the linear fractional map
\[\zeta \mapsto \frac{1+i\zeta}{1-i\zeta}\]
maps the (projectively) extended real line to the unit circle. Therefore, $\frac{1+i\zeta}{1-i\zeta}$ is a negative real number if and only if it is $-1$, and hence if  and only if $\zeta=\infty$. This concludes that $\varphi_\gamma$ is indeed analytic on $\mathbb{C}\setminus \mathbb{R}$.

Next, we view $\mathbb{R}$ as a contour oriented from $-\infty$ to $+\infty$. Then we have the following boundary values for $\varphi_\gamma$.
\begin{lemma}
\label{lem:behaviourphirealaxis}
We have that 
\begin{align} 
\varphi_{\gamma,\pm}(x) &= \pm \pi i\int_0^x d\mu_\gamma(s), &x\in(0,1), \label{eq:phiOn01} \\
\varphi_{\gamma,+}(x)&=\varphi_{\gamma,-}(x), &x<0, \label{eq:phinegative} 
\end{align}
whence $\varphi_\gamma$ can be analytically continued to form an analytic function $\varphi_\gamma: \mathbb{C}\setminus [0,\infty)\to\mathbb C$.
\end{lemma}
\begin{proof}
For \eqref{eq:phiOn01}, we note that we clearly have that $\varphi_{\gamma,-}(x)=-\varphi_{\gamma,+}(x)$ for all $x\in (0,1)$ by construction of $\varphi_\gamma$. Hence, we only prove the identity for $\varphi_{\gamma,+}$. Now, by considering the derivative and the value for $x=0$ of both sides of the following equation (using \eqref{eq:equilibriummeasuregamma}), we have that
\[\int_0^x d\mu_\gamma(s) = \sqrt{\frac{x}{\gamma}}
+\frac{2}{\pi}\arctan\left(\sqrt{\frac{\gamma-1}{\gamma}}\sqrt{\frac{x}{1-x}}\right) 
- \frac{2}{\pi}\sqrt{\frac{x}{\gamma}} \arctan\left(\sqrt{\frac{\gamma-1}{1-x}}\right), \qquad x\in (0,1).\]
Then using that $\arctan(y)+\arctan(1/y)=\frac{\pi}{2}$ if $y>0$, we obtain that
\begin{align*}
i\pi\int_0^x d\mu_\gamma(s) 
&=  2i\arctan\left(\sqrt{\frac{\gamma-1}{\gamma}}\sqrt{\frac{x}{1-x}}\right) 
+ 2i\sqrt{\frac{x}{\gamma}} \arctan\left(\sqrt{\frac{1-x}{\gamma-1}}\right)\\
&= \log\left(\frac{\sqrt\gamma\sqrt{1-x}+i\sqrt x\sqrt{\gamma-1}}{\sqrt\gamma\sqrt{1-x}-i\sqrt x\sqrt{\gamma-1}}\right)
+ \sqrt{\frac{x}{\gamma}} \log\left(\frac{\sqrt{\gamma-1}+i\sqrt{1-x}}{\sqrt{\gamma-1}-i\sqrt{1-x}}\right),
\end{align*}
where we used a standard identity between $\arctan$ and $\log$. The last expression is clearly equal to $\varphi_{\gamma,+}(x)$; this establishes \eqref{eq:phiOn01}. Since $\lim_{\varepsilon \rightarrow 0^+} \sqrt{x+i\varepsilon} = -\lim_{\varepsilon \rightarrow 0^+} \sqrt{x-i\varepsilon}$ for $x<0$, we conclude that \eqref{eq:phinegative} holds too.
\end{proof}

We now know that $\varphi_\gamma$ has the cut $[0,\infty)$. At the endpoint $z=0$ we have the following behaviour.
\begin{lemma} 
\label{lem:phionendpoints}
For the function $\varphi_\gamma$ defined by \eqref{eq:defphi}, we have the following:
\begin{equation*}
\varphi_\gamma(z) = \pm \frac{i}{\sqrt{c_\gamma}} \sqrt{z}+ O(z^{3/2}), \qquad \text{ as }z\to 0 \text{ and } \pm \operatorname{Im}(z)>0,
\end{equation*}
where $c_\gamma$ is as in \eqref{eq:defcgamma}. 
\end{lemma}
\begin{proof}
We only prove the behaviour of $\varphi_\gamma(z)$ around $z=0$ in the upper half plane; the behaviour on the lower half plane follows analogously. Hence, we are interested in the behaviour of
\begin{equation}
\label{eq:phiupperhalfplane}
\log\left(\frac{\sqrt\gamma\sqrt{1-z}+i\sqrt z\sqrt{\gamma-1}}{\sqrt\gamma\sqrt{1-z}-i\sqrt z\sqrt{\gamma-1}}\right)
+ \sqrt{\frac{z}{\gamma}} \log\left(\frac{\sqrt{\gamma-1}+i\sqrt{1-z}}{\sqrt{\gamma-1}-i\sqrt{1-z}}\right)
\end{equation}
around $z=0$. For this, we write
\[\zeta=\sqrt{\frac{\gamma-1}{\gamma}} \sqrt{\frac{z}{1-z}},\]
such that the first term of \eqref{eq:phiupperhalfplane} becomes $\log \frac{1+i\zeta}{1-i\zeta}$. We know that
\[\log \frac{1+i\zeta}{1-i\zeta} = 2i\zeta + O(\zeta^3), \qquad \zeta \rightarrow 0\]
so for the first term of \eqref{eq:phiupperhalfplane} we obtain
\[\log\left(\frac{\sqrt\gamma\sqrt{1-z}+i\sqrt z\sqrt{\gamma-1}}{\sqrt\gamma\sqrt{1-z}-i\sqrt z\sqrt{\gamma-1}}\right) = 2i \sqrt{\frac{\gamma -1}{\gamma}} \sqrt{z} + O(z^{3/2}), \qquad z\rightarrow 0.\]
For the second term of \eqref{eq:phiupperhalfplane} we have
\[\sqrt{\frac{z}{\gamma}} \log\left(\frac{\sqrt{\gamma-1}+i\sqrt{1-z}}{\sqrt{\gamma-1}-i\sqrt{1-z}}\right) = \sqrt{\frac{1}{\gamma}} \log\left(\frac{\sqrt{\gamma-1}+i}{\sqrt{\gamma-1}-i}\right) \sqrt{z} + O(z^{3/2}).\]
Hence, for $\operatorname{Im}(z)>0$ and $z\rightarrow 0$, we obtain
\begin{align*}
\varphi_\gamma(z) &= \left(2i\sqrt{\gamma-1} + \log\left(\frac{\sqrt{\gamma-1}+i}{\sqrt{\gamma-1} -i}\right)\right) \sqrt{\frac{z}{\gamma}} + O(z^{3/2}) \\
&=\frac{2i}{\sqrt{\gamma}} \left(\sqrt{\gamma-1} - \arctan(\sqrt{\gamma-1}) +\frac{\pi}{2}\right) \sqrt{z}+ O(z^{3/2}),
\end{align*}
which gives the desired result when invoking \eqref{eq:defcgamma}.
\end{proof}

\begin{figure}[h]
\centering
\resizebox{10.5cm}{6cm}{%
\begin{tikzpicture}[>=latex]
	\draw[-] (-7,0)--(7,0);
	\draw[-] (-4,-4)--(-4,4);
	\draw[fill] (-4,0) circle (0.1cm);
	\draw[fill] (3,0) circle (0.1cm);
	
	\node[above] at (-4.2,-0.075) {\large 0};	
	\node[above] at (2.5,0) {\large $1$};	
	
	\node[above] at (-0.2,2) {\large $\Delta_+$};
	\node[above] at (-0.2,-2.5) {\large $\Delta_-$};

	\draw[-,ultra thick] (-4,0)--(3,0);
	\draw[->, ultra thick] (-0.5,0) to (-0.3,0);

	\draw[-, ultra thick] (-4,0) to [out=60, in=120] (3,0);
	\draw[-, ultra thick] (-4,0) to [out=-60, in=-120] (3,0);

	\draw[->, ultra thick] (-0.5,1.77) to (-0.3,1.77);
	\draw[->, ultra thick] (-0.5,-1.77) to (-0.3,-1.77);
	
	\draw[-, color=white] (-4,3)--(-4,4);
	\draw[-, color=white] (-4,-4)--(-4,-3);
\end{tikzpicture}
}
\caption{The set of contours $\Sigma_S$.}
\label{fig:contoursS}
\end{figure}
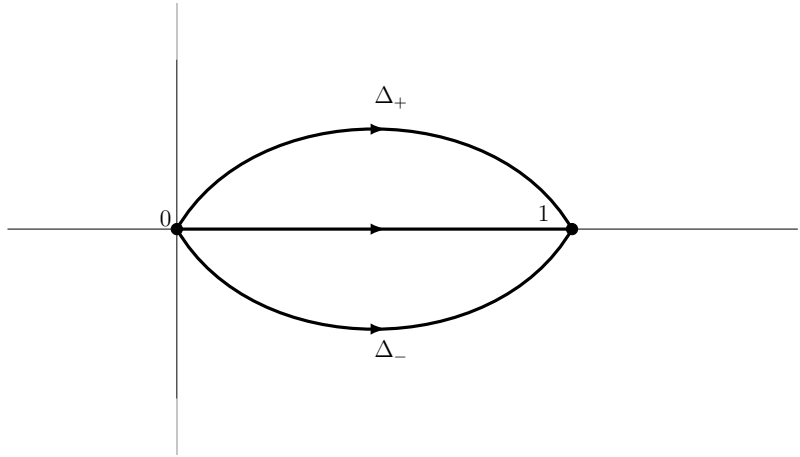

Next we open a lens from $0$ to $1$. This means that we take two contours $\Delta_+$ and $\Delta_-$, both going from $0$ to $1$, where $\Delta_+$ goes through the upper half plane and $\Delta_-$ through the lower half plane. We write $\Sigma_S=(0,1) \cup \Delta_+ \cup \Delta_-$ for the collection of contours that we now have under consideration, see Figure \ref{fig:contoursS}.

As is customary, we refer to the region enclosed between $\Delta_+$ and $\Delta_-$ as the interior of the lens, the region between $\Delta_+$ and $(0,1)$ as the upper part of the lens and between $(0,1)$ and $\Delta_-$ as the lower part of the lens. We then define the function $S: \mathbb{C}\setminus \Sigma_S \rightarrow \mathbb{C}^{2\times 2}$ by
\begin{equation}
\label{eq:defS}
S(z) = \left\{
\begin{array}{ll}  T(z) \begin{pmatrix} 1 & 0\\ - z^{-\nu} e^{2n\varphi_\gamma(z)} & 1\end{pmatrix} & \text{ in the upper part of the lens},\\
T(z) \begin{pmatrix} 1 & 0\\ z^{-\nu} e^{2n\varphi_\gamma(z)} & 1\end{pmatrix} & \text{ in the lower part of the lens}, \\
T(z) & \text{ outside the lens}. \end{array}
\right.
\end{equation}
Then $S$ satisfies the following Riemann-Hilbert problem.
\begin{itemize}
\item[RH-S1] $S:\mathbb{C}\setminus \Sigma_S \to \mathbb{C}^{2\times 2}$ is analytic.
\item[RH-S2a]On $(0,1)$ we have the following jump:
\begin{equation*}
S_+(x) = S_-(x) \begin{pmatrix} 0 & x^\nu \\ -x^{-\nu} & 0\end{pmatrix} \qquad x\in (0,1).
\end{equation*}
\item[RH-S2b] On the lips $\Delta_+$ and $\Delta_-$ we have the following jump, where $S_+$ and $S_-$ are determined by the orientation of $\Delta_+$ and $\Delta_-$:
\begin{equation*}
S_+(z) = S_-(z) \begin{pmatrix} 1 & 0\\ z^{-\nu} e^{2n\varphi_\gamma(z)} & 1\end{pmatrix} \qquad z\in \Delta_+ \cup \Delta_-.
\end{equation*}
\item[RH-S3] As $z\to\infty$ we have $S(z) = I + O(1/z)$. 
\item[RH-S4a] As $z\to 0$ we have
\begin{align}
S(z) &= O \begin{pmatrix} 1 & h_{\nu}(z)\\ 1 & h_{\nu}(z)\end{pmatrix} & \text{ for $z$ outside the lens.}\\
S(z) &= O \begin{pmatrix} h_{-\nu}(z) & h_\nu(z)\\ h_{-\nu}(z) & h_\nu(z)\end{pmatrix} & \text{ for $z$ inside the lens.}
\end{align}
\item[RH-S4b] As $z\to 1$ we have
\begin{align}
S(z) &= O \begin{pmatrix} 1 & \log(z-1)\\ 1 & \log(z-1)\end{pmatrix} & \text{ for $z$ outside the lens.}\\
S(z) &= O \begin{pmatrix} \log(z-1) & \log(z-1)\\ \log(z-1) & \log(z-1)\end{pmatrix} & \text{ for $z$ inside the lens.} \label{eq:S1inside}
\end{align}
\end{itemize}
The only non-trivial properties to prove is the behaviour inside the lens in RH-S4a and RH-S4b. For the first, note that $z^{-\nu}h_\nu(z) = h_{-\nu}(z)$ by definition of $h_\nu$ \eqref{eq:defhnu}. Furthermore, by Lemma \ref{lem:phionendpoints}, we have that $e^{2 n \varphi_\gamma(z)} = O(1)$ as $z\to 0$, whence $O(1\pm h_\nu(z) e^{2n\varphi_\gamma(z)})=O(h_\nu(z))$, regardless of the exact value of $\nu$. So, for $z$ inside the lens we have as $z\to 0$
\begin{align*}
S(z) &= O\begin{pmatrix} 1 & h_{\nu(z)}\\ 1 & h_{\nu}(z)\end{pmatrix} 
\begin{pmatrix} 1 & 0\\ \mp z^{-\nu} e^{2n\varphi_\gamma(z)} & 1\end{pmatrix} 
= O\begin{pmatrix} 1\mp h_{-\nu}(z) e^{2 n\varphi_\gamma(z)} & h_\nu(z)\\ 1\mp h_{-\nu}(z) e^{2 n\varphi_\gamma(z)} & h_\nu(z)\end{pmatrix} \\
&= O \begin{pmatrix} h_{-\nu}(z) & h_\nu(z)\\ h_{-\nu}(z) & h_\nu(z)\end{pmatrix}.
\end{align*}
For the behaviour inside the lens as $z\rightarrow 1$ in RH-S4b, we remark that $z^{-\nu} e^{2 n\varphi_\gamma(z)}\to 1$ as $z\to 1$. Combining this with RH-T4b and the definition of S \eqref{eq:defS}, yields \eqref{eq:S1inside}.

\subsection{Global parametrix}
For large $n$ the jump matrices of $S$ on the lips of the lens are close to the unit matrix, see Lemma \ref{lem:realpartoncontours}. The global parametrix problem for the Riemann-Hilbert problem neglects these jumps altogether and is hence the following:
\begin{itemize}
\item[RH-N1] $N:\mathbb{C}\setminus [0,1]\to\mathbb{C}^{2\times 2}$ is analytic.
\item[RH-N2] We have the following jump for $x\in (0,1)$ (oriented away from the origin)
\[N_+(x) = N_-(x) \begin{pmatrix} 0 & x^\nu\\ - x^{-\nu} & 0\end{pmatrix}\]
\item[RH-N3] As $z\to\infty$ we have $N(z) = I+O(1/z)$.
\end{itemize}

Since we aim to approximate $S$ by $N$ sufficiently far away from the end points, we leave some freedom for the behavior of $N$ around $z=0$ and $z=1$. A solution is readily available in the literature \cite[Chapter 5]{Kuijlaars_2004}. After an appropriate translation it yields
\begin{align}
\label{eq:defN}
N(z) &= 2^{-\nu\sigma_3}
\begin{pmatrix}
\displaystyle\frac{a(z)+a(z)^{-1}}{2} & \displaystyle\frac{a(z)-a(z)^{-1}}{2i}\\
\displaystyle\frac{a(z)-a(z)^{-1}}{-2i} & \displaystyle\frac{a(z)+a(z)^{-1}}{2}
\end{pmatrix}
\left(1+\sqrt{\frac{z-1}{z}}\right)^{\nu\sigma_3},
\end{align}
where $a(z) = \left(\frac{z-1}{z}\right)^\frac{1}{4}$. Here and later, we adhere to the usual notation for the third Pauli matrix
\[\sigma_3=\begin{pmatrix}
1 & 0 \\
0 & -1 \\
\end{pmatrix}.\]

With this definition of $N$ we should add the following local behaviors.
\begin{itemize}
\item[RH-N4a] As $z\to 0$ we have
\[N(z) = z^{-1/4} O \begin{pmatrix} z^{-\nu/2} & z^{\nu/2}\\ z^{-\nu/2} & z^{\nu/2}\end{pmatrix}.\]
\item[RH-N4b] As $z\to 1$ we have
\[N(z) = O ((z-1)^{-1/4}).\]
\end{itemize}

\subsection{Local parametrices}
We have two local parametrix problems; one around $z=0$ and one around $z=1$. For finding the solutions, we rely on the results available in \cite{Kuijlaars_2004}, where a similar local parametrix problem was studied. We only explicitly do this analysis for the local parametrix around $z=0$, since the other local parametrix problem is similar and its details are not needed for our further analysis.

\subsubsection{The local parametrix around the origin}
The local parametrix problem around $z=0$ is the following Riemann-Hilbert problem, where $r>0$ is a small number.
\begin{itemize}
\item[RH-P1] $P: \overline{D(0,r)}\setminus \Sigma_S \rightarrow \mathbb{C}^{2\times 2}$ is analytic.
\item[RH-P2a] On $(0,r)$ we have the following jump:
\begin{equation*}
P_+(x) = P_-(x) \begin{pmatrix} 0 & x^\nu \\ -x^{-\nu} & 0\end{pmatrix}, \qquad x\in (0,r).
\end{equation*}
\item[RH-P2b] On the contours $\Delta_+\cap D(0,r)$ and $\Delta_-\cap D(0,r)$ we have the following jump:
\begin{equation*}
P_+(z) = P_-(z) \begin{pmatrix} 1 & 0\\ z^{-\nu} e^{2n\varphi_\gamma(z)} & 1\end{pmatrix}, \qquad z\in (\Delta_+ \cup \Delta_-)\cap D(0,r).
\end{equation*}
\item[RH-P3] As $z\to 0$ we have
\begin{align}
P(z) &= O \begin{pmatrix} 1 & h_{\nu}(z)\\ 1 & h_{\nu}(z)\end{pmatrix}, & \text{ for $z$ outside the lens,}\\
P(z) &= O \begin{pmatrix} h_{-\nu}(z) & h_\nu(z)\\ h_{-\nu}(z) & h_\nu(z)\end{pmatrix}, & \text{ for $z$ inside the lens.}
\end{align}
\item[RH-P4] We have that $P(z)=(I+O(1/n))N(z)$ uniformly for $\lvert z\rvert = r$ as $n\rightarrow \infty$.
\end{itemize}
To find the solution to this problem, we make use of the results in \cite{Kuijlaars_2004}, where a similar local parametrix problem was solved. For this, we consider the following two rays in the complex plane:
\begin{align}
\eta_+ &= \{t e^{\frac{\pi i}{3}} \mid t> 0\},\\
\eta_- &= \{t e^{\frac{-\pi i}{3}} \mid t> 0\}.
\end{align}
We see these rays as being oriented away from the origin. Then, we consider the function $B: \mathbb{C}\setminus ([0,\infty) \cup \eta_+ \cup \eta_-) \rightarrow \mathbb{C}^{2\times 2}$ defined by
\begin{align}
\label{eq:defB}
B(z) &=\left\{\begin{array}{ll}
\begin{pmatrix}
\frac{1}{2} H_\nu^{(2)}\left(2\sqrt z\right) & \frac{1}{2} H_\nu^{(1)}\left(2\sqrt z\right)\\
\pi \sqrt{-z} (H_\nu^{(2)})'(2 \sqrt z) & \pi \sqrt{-z} (H_\nu^{(1)})'(2\sqrt z) 
\end{pmatrix} e^{-\frac{\pi i\nu}{2} \sigma_3}, & 0 < \arg z < \frac{\pi}{3},\\
\quad\\
\begin{pmatrix}
I_\nu(2\sqrt{-z}) & -\frac{i}{\pi} K_\nu(2\sqrt{-z})\\
-2\pi i \sqrt{-z} I_\nu'(2\sqrt{-z}) & -2\sqrt{-z} K_\nu'(2\sqrt{-z})
\end{pmatrix}, & |\arg z|>\frac{\pi}{3},\\
\quad\\
\begin{pmatrix}
\frac{1}{2} H_\nu^{(1)}\left(2\sqrt z\right) & -\frac{1}{2} H_\nu^{(2)}\left(2\sqrt z\right)\\
-\pi \sqrt{-z} (H_\nu^{(1)})'(2 \sqrt z) & \pi \sqrt{-z} (H_\nu^{(2)})'(2\sqrt z) 
\end{pmatrix} e^{\frac{\pi i\nu}{2} \sigma_3}, & -\frac{\pi}{3} < \arg z < 0.
\end{array}\right.
\end{align}
Here $I_\nu$ and $K_\nu$ are modified Bessel functions and $H_\nu^{(1)}$ and $H_\nu^{(2)}$ are Hankel functions. See \cite[Chapter 9]{Abramowitz_Stegun} for more details on these functions. Indeed, if we compare \eqref{eq:defB} with the definitions on page 367 of \cite{Kuijlaars_2004}, we see that their function $\Psi$ is related to our function $B$ by 
\begin{equation}
\label{eq:relationBPsi}
B(z)=\sigma_3 \Psi(-z) \sigma_3.
\end{equation}

Using \eqref{eq:relationBPsi} and the results of \cite{Kuijlaars_2004}, we directly obtain that $B$ satisfies the following Riemann-Hilbert problem:
\begin{itemize}
\item[RH-B1] $B: \mathbb{C}\setminus ([0,\infty) \cup \eta_+ \cup \eta_-) \rightarrow \mathbb{C}^{2\times 2}$ is analytic.
\item[RH-B2] $B$ has the following jumps (all contours oriented away from the origin)
\begin{align}
B_+(x) &= B_-(x) \begin{pmatrix} 0 & 1\\ -1 & 0\end{pmatrix} & x\in (0,\infty)\\
B_+(z) &= B_-(z) \begin{pmatrix} 1 & 0\\ e^{\mp \pi i\nu} & 1\end{pmatrix} & z\in \eta_\pm.
\end{align}
\item[RH-B3] As $z\to 0$ we have
\begin{align}
B(z) &= z^{\nu/2} O \begin{pmatrix} 1 & h_{-\nu(z)}\\ 1 & h_{-\nu}(z) \end{pmatrix} & \lvert \arg(z)-\pi \rvert < \frac{2\pi}{3}\\
B(z) &= z^{\nu/2} O \begin{pmatrix} h_{-\nu}(z) & h_{-\nu}(z)\\ h_{-\nu}(z) & h_{-\nu}(z) \end{pmatrix} & 0<\lvert\arg(z)\rvert<\frac{\pi}{3}.
\end{align}
\end{itemize}
We note that $B$ is not the unique solution of the problem; requiring a certain asymptotic behaviour as $z\rightarrow \infty$ would make it unique. For the particular definition of $B$ in \eqref{eq:defB}, by \eqref{eq:relationBPsi} and the results in \cite{Kuijlaars_2004}, we have that 
\begin{equation}
\label{eq:Binfty}
B(z) = \left(2\pi \sqrt{-z}\right)^{-\frac{\sigma_3}{2}} 
\left(\frac{1}{\sqrt 2} \begin{pmatrix} 1 & -i\\ -i & 1\end{pmatrix} + \mathcal O\left(\frac{1}{\sqrt z}\right)\right) e^{2\sqrt{-z}\sigma_3}, \qquad z\rightarrow \infty.
\end{equation}
This function $B$ is not directly the solution that we need in our analysis; we need to transform $B$. For this, we write $D(0,r)$ and $\overline{D(0,r)}$ for the open and closed disc of radius $r$ around $z=0$, respectively. Then we define the function $f_\gamma: D(0,1) \rightarrow \mathbb{C}$ by
\begin{equation}
\label{eq:deffgamma}
f_\gamma(z) = -\frac{1}{4} \varphi_\gamma(z)^2,
\end{equation} 
where $\varphi_\gamma$ is as in \eqref{eq:defphi}. Although $\varphi_\gamma$ has $[0,\infty)$ as a cut, $f_\gamma$ is analytic on $D(0,1)$ due to the square in its definition and the jump \eqref{eq:phiOn01} of $\varphi_\gamma$. Then, by Lemma \ref{lem:phionendpoints} and \eqref{eq:defcgamma}, we have that
\begin{equation}
\label{eq:fgamma0}
f_\gamma(z) = \frac{\pi^2}{4 c_\gamma}  z + \mathcal O(z^2), \qquad \text{ as } z\rightarrow 0.
\end{equation}
In particular, there exists an $0<r_0<\frac{1}{2}$ such that the derivative of $f_\gamma$ does not vanish on $D(0,r_0)$. Then the same holds for every $0<r\leq r_0$ and hence the restriction of $f_\gamma$ to $D(0,r)$ is a conformal map for every such $r$. In what follows, we assume that $0<r\leq r_0$.

We use the definition \eqref{eq:defN} of the matrix function $N$ and the above definition \eqref{eq:deffgamma} of $f_\gamma$, to define the function
\begin{equation}
\label{eq:defEn}
E_n(z) = N(z) 
(-z)^{\frac{\nu}{2}\sigma_3}
\frac{1}{\sqrt 2} \begin{pmatrix} 1 & i\\ i & 1\end{pmatrix} (2\pi n)^\frac{\sigma_3}{2} \left(-f_\gamma(z)\right)^{\frac{\sigma_3}{4}},
\end{equation}
which we want to consider for $z\in D(0,r)$. In \cite{Kuijlaars_2004} a similar function is defined which is also denoted by $E_n$. Analogous to that, we have the following.
\begin{lemma}
The function $E_n$ defined by \eqref{eq:defEn} defines an analytic and non-singular function on the disc $D(0,r)$. 
\end{lemma}

We remark that when we first used the contours $\Delta_+$ and $\Delta_-$ in subsection \ref{subsec:openinglens}, we only required that $\Delta_+$ went from $0$ to $1$ through the upper half plane, and $\Delta_-$ similarly through the lower half plane. Since our previous results hold for any such $\Delta_+$ and $\Delta_-$, we may assume more about these contours. Namely, since the function $f_\gamma$ defined in \eqref{eq:deffgamma} is conformal, and maps positive numbers to positive numbers and $0$ to $0$, we may assume that
\begin{align}
f_\gamma(\Delta_+ \cap D(0,r)) &\subset \eta_+, \label{eq:Delta+eta+} \\
f_\gamma(\Delta_- \cap D(0,r)) &\subset \eta_- \label{eq:Delta-eta-},
\end{align}  
i.e. $f_\gamma$ maps the contours $\Sigma_S \cap D(0,r)$ to the cuts of the function $B$. Hence for any $n\geq 1$ we can study the function $Q: \overline{D(0,r)} \setminus \Sigma_S$ defined by 
\begin{equation}
\label{eq:defQ}
Q(z)= E_n(z) B(n^2 f_\gamma(z)), \qquad z\in \overline{D(0,r)} \setminus \Sigma_S,
\end{equation}
where again we suppress the dependency on $n$ from the notation. It immediately follows that $Q$ satisfies the following Riemann-Hilbert problem:
\begin{itemize}
\item[RH-Q1] $Q : \overline{D(0,r)}\setminus \Sigma_S \to\mathbb{C}^{2\times 2}$ is analytic.
\item[RH-Q2] $Q$ has the following jumps (all contours oriented away from the origin)
\begin{align}
Q_+(x) &= Q_-(x) \begin{pmatrix} 0 & 1\\ -1 & 0\end{pmatrix}, & x\in (0,r),\\
Q_+(z) &= Q_-(z) \begin{pmatrix} 1 & 0\\ e^{\mp \pi i\nu} & 1\end{pmatrix}, & z\in \Delta^\pm \cap D(0,r).
\end{align}
\item[RH-Q3] As $z\to 0$ we have 
\begin{align}
Q(z) &= z^{\nu/2} O \begin{pmatrix} 1 & h_{-\nu(z)}\\ 1 & h_{-\nu}(z) \end{pmatrix}, & \text{ outside the lens,}\\
Q(z) &= z^{\nu/2} O \begin{pmatrix} h_{-\nu}(z) & h_{-\nu}(z)\\ h_{-\nu}(z) & h_{-\nu}(z) \end{pmatrix}, & \text{ inside the lens.}
\end{align}
\item[RH-Q4] We have
\[Q(z) = (I+O(1/n)) N(z)
(-z)^{\frac{\nu}{2}\sigma_3} e^{-n\varphi(z) \sigma_3}.\] 
uniformly on $|z|=r$ as $n\to\infty$.
\end{itemize}
%
%
%
Using this, one straightforwardly checks that the function $P: \overline{D(0,r)}\setminus \Sigma_S \rightarrow \mathbb{C}^{2\times 2}$, defined by
\begin{equation}
\label{eq:defP}
P(z) = Q(z) e^{n\varphi_\gamma(z)\sigma_3} (-z)^{-\frac{\nu}{2}\sigma_3}, \qquad z\in \overline{D(0,r)}\setminus \Sigma_S,
\end{equation}
satisfies the local parametrix problem RH-P.
%
A last property of $P$ that we need in the further analysis is the following fact.
\begin{lemma}
\label{lem:SP-1}
We have that $S(z)P(z)^{-1}$ is analytic in $D(0,r)$.
\end{lemma}
\begin{proof}
First of all, since the jumps of $P$ that appear in RH-P2a and RH-P2b are the same as the jumps of $S$ that appear in RH-S2a and RH-S2b, we conclude that the only potentially singular point is $z=0$. Then, using the property RH-B3 of $B$, we have that as $z\to 0$ outside the lens
\begin{multline}
P(z)^{-1} = (-z)^{\frac{\nu}{2}\sigma_3} e^{-n\varphi_\gamma(z)\sigma_3}
B\left(n^2 f_\gamma(z)\right)^{-1} E_n(z)^{-1}\\
= \begin{pmatrix} z^{\nu/2} & 0\\ 0 & z^{-\nu/2} \end{pmatrix}
z^{\nu/2} O\begin{pmatrix} h_{-\nu}(z) & h_{-\nu}(z)\\ 1 & 1\end{pmatrix}
= O\begin{pmatrix} h_{\nu}(z) & h_{\nu}(z)\\ 1 & 1\end{pmatrix}.
\end{multline}
Then using RH-S4a we have outside the lens as $z\to 0$ that 
\begin{align}
S(z) P(z)^{-1} &= O\left(
\begin{pmatrix} 1 & h_\nu(z)\\ 1 & h_\nu(z)\end{pmatrix}
\begin{pmatrix} h_{\nu}(z) & h_{\nu}(z)\\ 1 & 1\end{pmatrix}\right) = O(h_\nu(z)).
\end{align}
By the definition \eqref{eq:defhnu} of $h_\nu$, this implies that $S(z)P(z)^{-1}$ has no negative powers in its Laurent series around $z=0$. We conclude that $S(z) P(z)^{-1}$ is analytic in $z=0$.
\end{proof}

\subsubsection{The second local parametrix}
The local parametrix problem around $z=1$ can be solved similarly as the local parametrix problem around $z=0$. More specifically, if we use the function $\Psi$ from \cite{Kuijlaars_2004} for the parameter $\alpha=0$ and transform it as we did in the previous section, we have the solution to the local parametrix problem. Namely, we find that that there is a $0<r_1<\frac{1}{2}$ such that for all $0<r\leq r_1$, we may assume that the contours $\Delta_+$ and $\Delta_-$ are chosen such that there exists a function $\tilde{P}: \overline{D(1,r)}\setminus \Sigma_S \rightarrow \mathbb{C}^{2\times 2}$, that satisfies the following Riemann-Hilbert problem.
\begin{itemize}
\item[RH-$\tilde P$1] $\tilde{P} : \overline{D(1,r)}\setminus \Sigma_S\to \mathbb{C}^{2\times 2}$ is analytic.
\item[RH-$\tilde{P}$2a] On $(1-r,1)$ we have the following jump:
\begin{equation*}
\tilde{P}_+(x) = \tilde{P}_-(x) \begin{pmatrix} 0 & x^\nu \\ -x^{-\nu} & 0\end{pmatrix}, \qquad x\in (0,r).
\end{equation*}
\item[RH-P2b] On the contours $\Delta_+\cap D(1,r)$ and $\Delta_-\cap D(1,r)$ we have the following jump:
\begin{equation*}
\tilde{P}_+(z) = \tilde{P}_-(z) \begin{pmatrix} 1 & 0\\ z^{-\nu} e^{2n\varphi_\gamma(z)} & 1\end{pmatrix}, \qquad z\in (\Delta_+ \cup \Delta_-)\cap D(1,r).
\end{equation*}
\item[RH-$\tilde P$3] As $z\rightarrow 1$, we have 
\begin{align}
\tilde{P}(z) &= O \begin{pmatrix} 1 & \log(z-1)\\ 1 & \log(z-1)\end{pmatrix}, & \text{ for $z$ outside the lens.}\\
\tilde{P}(z) &= O \begin{pmatrix} \log(z-1) & \log(z-1))\\ \log(z-1) & \log(z-1)\end{pmatrix}, & \text{ for $z$ inside the lens.} 
\end{align}
\item[RH-$\tilde P$4] Matching: $\tilde P(z) = (I+O(1/n)) N(z)$ uniformly for $|z-1|=r$ as $n\to\infty$, where $N$ is defined in \eqref{eq:defN}.
\end{itemize}

Since the computations to arrive at this solution are standard by the work of \cite{Kuijlaars_2004} and we do not need the explicit form of this solution in our further analysis, we omit further details. However, we do note that we have the following, which is analogous to Lemma \ref{lem:SP-1}.
\begin{lemma}
\label{lem:StildeP-1}
We have that $S(z)\tilde{P}(z)^{-1}$ is analytic in $D(1,r)$.
\end{lemma}

\subsection{Final transformation}
\label{subsec:finaltransformation}
We can now make the final transformation: informally, we show that if one `glues' the global parametrix solution $N$ to the local solutions $P$ and $\tilde{P}$, one asymptotically obtains the solution for the Riemann-Hilbert problem RH-S. We have the following result.
\begin{lemma}
\label{lem:realpartoncontours}
There exists an $r<\min\{r_0,r_1\}$ and an opening of the lips of the lens such that
\begin{align} \label{eq:realpartoncontours}
\operatorname{Re}(\varphi_\gamma(z))&\leq-c, & z \in (\Delta_+ \cup \Delta_-) \setminus (D(0,r)\cup D(1,r)),
\end{align}
for some constant $c>0$, where $\varphi_\gamma$ is defined by \eqref{eq:defphi}.
\end{lemma}

\begin{proof}
By $\eqref{eq:phiOn01}$, we have that $\varphi_{\gamma,\pm}$ is a purely imaginary function on $(0,1)$. In fact, $\operatorname{Im}(\varphi_{\gamma,+})$ is strictly increasing, and $\operatorname{Im}(\varphi_{\gamma,-})$ is strictly decreasing on $(0,1)$. By the Cauchy-Riemann equations $\operatorname{Re}(\varphi_\gamma(z))$ will be negative on any subset of the strip $0<\operatorname{Re}(z)<1$ that is close enough to $(0,1)$, but does not contain any point of $(0,1)$. Combining this with the continuity of $\varphi_\gamma$ yields the constant $c>0$.
\end{proof}

\begin{figure}[h]
\centering
\resizebox{10.5cm}{6cm}{%
\begin{tikzpicture}[>=latex]
	\draw[-] (-7,0)--(7,0);
	\draw[-] (-4,-4)--(-4,4);
	\draw[fill] (-4,0) circle (0.1cm);
	\draw[fill] (3,0) circle (0.1cm);
	\draw[ultra thick] (-4,0) circle (1.5cm);
	\draw[ultra thick] (3,0) circle (1.5cm);
	\node[above] at (-4.2,0) {\large 0};	
	\node[above] at (-2.2,0) {\large $r$};
	\node[above] at (3.3,0) {\large $1$};	
	\node[above] at (.95,0) {\large $1-r$};	
	
	\node[above] at (0.2,1.8) {\large $\Delta_+$};
	\node[above] at (0.2,-2.5) {\large $\Delta_-$};

	\draw[-,ultra thick] (-4,0)--(3,0);
	\draw[-, ultra thick] (-4,0) to [out=60, in=120] (3,0);
	\draw[-, ultra thick] (-4,0) to [out=-60, in=-120] (3,0);

	\draw[->, ultra thick] (-5.06,1.08) to (-5.16,0.96);
	\draw[->, ultra thick] (4.14,0.96) to (4.04,1.08);
	
	\draw[->, ultra thick] (-0.5,1.78) to (-0.3,1.78);
	\draw[->, ultra thick] (-0.5,-1.78) to (-0.3,-1.78);
	
	\draw[-, color=white] (-4,3)--(-4,4);
	\draw[-, color=white] (-4,-4)--(-4,-3);
\end{tikzpicture}
}
\caption{The set of contours $\Sigma_R$.}
\label{fig:contoursR}
\end{figure}
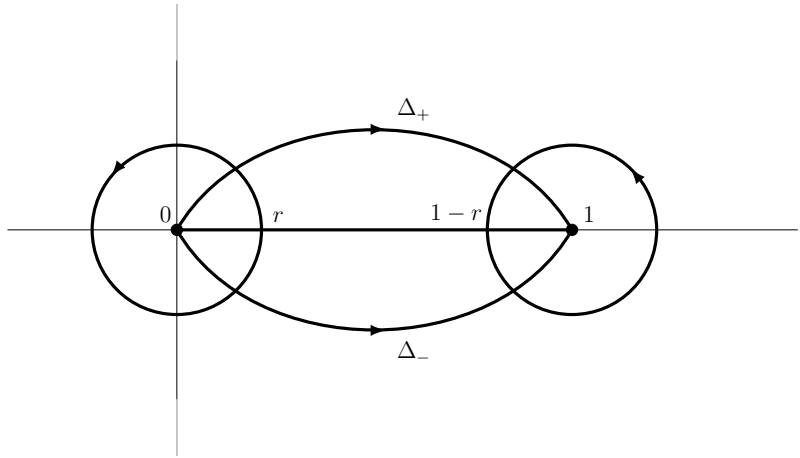

Now let us fix such an $r$ and opening. We collect all the cuts that we have used in the set of contours
\[\Sigma_R=(0,1) \cup \Delta_+ \cup \Delta_- \cup \partial D(0,r) \cup \partial D(1,r),\]
where $\partial D(0,r)$ and $\partial D(1,r)$ are the boundaries of the disks $D(0,r)$ and $D(1,r)$, respectively, see Figure \ref{fig:contoursR}.

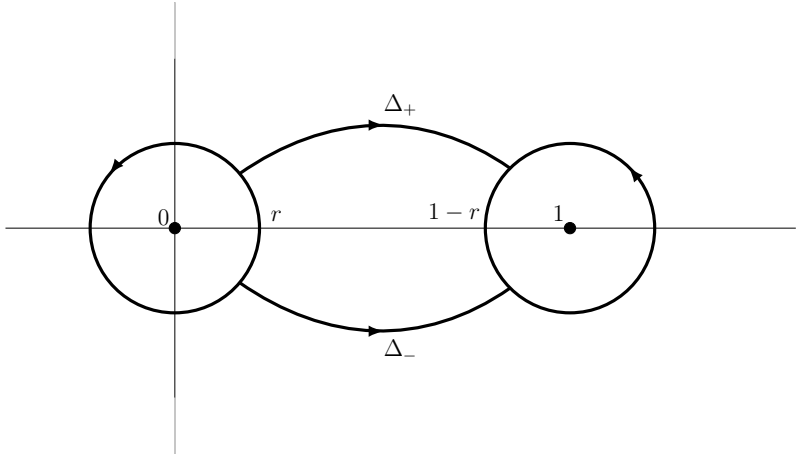
\begin{figure}[h]
\centering
\resizebox{10.5cm}{6cm}{%
\begin{tikzpicture}[>=latex]
	\draw[-] (-7,0)--(7,0);
	\draw[-] (-4,-4)--(-4,4);
	\draw[fill] (-4,0) circle (0.1cm);
	\draw[fill] (3,0) circle (0.1cm);
	\draw[ultra thick] (-4,0) circle (1.5cm);
	\draw[ultra thick] (3,0) circle (1.5cm);
	\node[above] at (-4.2,-0.075) {\large 0};	
	\node[above] at (-2.2,0) {\large $r$};
	\node[above] at (2.8,0) {\large $1$};	
	\node[above] at (.95,0) {\large $1-r$};	
	
	\node[above] at (0,1.85) {\large $\Delta_+$};
	\node[above] at (0,-2.5) {\large $\Delta_-$};

	\draw[-, ultra thick] (-2.86,0.96) to [out=35, in=145] (1.939,1.061);
	\draw[-, ultra thick] (-2.86,-0.96) to [out=-35, in=-145] (1.939,-1.061);

	\draw[->, ultra thick] (-5.06,1.08) to (-5.16,0.96);
	\draw[->, ultra thick] (4.14,0.96) to (4.04,1.08);
	
	\draw[->, ultra thick] (-0.5,1.825) to (-0.3,1.825);
	\draw[->, ultra thick] (-0.5,-1.825) to (-0.3,-1.825);
	
	\draw[-, color=white] (-4,3)--(-4,4);
	\draw[-, color=white] (-4,-4)--(-4,-3);
\end{tikzpicture}
}
\caption{The set of contours $\Sigma'_R$.}
\label{fig:contoursRprime}
\end{figure}

Then we define the function $R: \mathbb{C}\setminus \Sigma_R \rightarrow \mathbb{C}^{2\times 2}$, by
\begin{equation}
\label{eq:defR}
R(z) = \left\{\begin{array}{ll}
S(z) P(z)^{-1} & |z|<r\\
S(z) \tilde P(z)^{-1} & |z-1|<r\\
S(z) N(z)^{-1} & \text{otherwise.}\end{array}\right.
\end{equation}
We immediately note that the jumps of $R$ inside the disk $D(0,r)$ disappear due to the combination of RH-S2 and RH-P2. Likewise, the jumps inside the disk $D(1,r)$ disappear due to RH-S2 and RH-$\tilde{P}$2. Furthermore, by Lemma \ref{lem:SP-1} and Lemma \ref{lem:StildeP-1}, $R$ is analytic in the (potentially singular) points $z=0$ and $z=1$. Lastly, we also note that by combining RH-N2 with RH-S2a, we see that the (potential) jump of $R$ across $(r,1-r)$ does not exist. Therefore, we see $R$ as being defined on $\mathbb{C}\setminus \Sigma'_R$, where
\[\Sigma'_R=(\Delta_+ \cup \Delta_- \cup \partial D(0,r) \cup \partial D(1,r)) \setminus(D(0,r) \cup D(1,r)) \cup \partial D(1,r)\]
is the set of remaining contours, see Figure \ref{fig:contoursRprime}.

Combining the definition of $R$ with the Riemann-Hilbert problems of $S$, $N$, $P$ and $\tilde{P}$, we then immediately have that $R$ satisfies the following Riemann-Hilbert problem.
\begin{itemize}
\item[RH-R1] $R: \mathbb{C}\setminus \Sigma'_R \rightarrow \mathbb{C}^{2\times 2}$ is analytic.
\item[RH-R2a] For $z\in \partial D(0,r) \cup \partial D(1,r)$, we have that
\[R_+(z) = R_-(z) (I + O(1/n)).\]
\item[RH-R2b] For $z\in (\Delta_+ \cup \Delta_-) \setminus(D(0,r) \cup D(1,r))$, we have that
\begin{equation*}
R_+(z) = R_-(z) \begin{pmatrix} 1 & 0\\ z^{-\nu} e^{2n\varphi_\gamma(z)} & 1\end{pmatrix}.
\end{equation*}
\item[RH-R3] As $z\to\infty$ we have $R(z) = I+O(1/z)$.
\end{itemize}

Therefore, by \eqref{eq:realpartoncontours}, we have the following asymptotic behaviour for the jump matrix in RH-R2b:
\[\begin{pmatrix} 1 & 0\\ z^{-\nu} e^{2n\varphi_\gamma(z)} & 1\end{pmatrix} = I+O\left(e^{-2nc}\right), \qquad n\to \infty.\]
Since the jump matrices on $\Sigma'_R$ tend to the identity as $n\rightarrow \infty$, we can impose the general result concerning Riemann-Hilbert problems (see for example \cite{Bleher_Liechty,Deift}), to conclude that
\begin{equation}
\label{eq:asymptoticR}
R(z)=I+O(1/n), \qquad n\rightarrow \infty,
\end{equation}
uniformly for $z\in \mathbb{C}\setminus \Sigma'_R$. In the following, we use this asymptotic result to obtain an asymptotic result for the original matrix $Y$, that we need to obtain the asymptotics for the kernel that has our interest, according to \eqref{eq:kernelinY}.

\subsection{Asymptotics of the correlation kernel}
We have all the ingredients necessary to complete the proof of Proposition \ref{prop:resultRHprob}. As indicated in Section \ref{sec:outline}, this is the last thing that we need to prove.

\begin{proof}[Proof of Proposition \ref{prop:resultRHprob}]
If we invert the steps $Y \mapsto T \mapsto S \mapsto R$ in our Riemann-Hilbert problem, we obtain that for $z$ in the upper lens such that $\lvert z \rvert <r$, we have that
\begin{equation}
\label{eq:YafterQ}
\begin{pmatrix} Y_{11}(z) \\ Y_{21}(z)\end{pmatrix} =z^{-\frac{\nu}{2}} e^{n(g_\gamma(z)-\frac{\ell_\gamma}{2}+\varphi_\gamma(z))} e^{\frac{n\ell_\gamma}{2}\sigma_3} R(z) E_n(z) B(n^2 f_\gamma(z)) 
e^{\frac{\pi i\nu}{2} \sigma_3} \begin{pmatrix} 1 \\ 1\end{pmatrix}.
\end{equation}
Notice that, due to the definition \eqref{eq:defB} of $B$ for $0 < \arg z < \frac{\pi}{3}$, and the well-known identity \cite[Eq. 10.4.4]{DLMF} that connects Hankel and Bessel functions
\[H_{\nu}^{(1)}(z) + H_{\nu}^{(2)}(z)=2J_\nu(z),\]
we have for $x\in (0,r)$ that
\[B_+(n^2f_\gamma(x))
e^{\frac{\pi i\nu}{2} \sigma_3} \begin{pmatrix} 1 \\ 1\end{pmatrix} = \begin{pmatrix}
J_\nu(2n\sqrt{f_\gamma(x)}) \\
2\pi n \sqrt{-f_\gamma(x)} J'_\nu(2n\sqrt{f_\gamma(x)})
\end{pmatrix}.\]
We note that by \eqref{eq:gplus+minusinterval}, \eqref{eq:gplus-minusinterval} and \eqref{eq:phiOn01} we have for $x\in (0,1)$ that
\[g_{\gamma,+}(x)-\frac{\ell_\gamma}{2}+\varphi_{\gamma,+}(x)=\frac{V_\gamma(x)}{2} + \pi i,\]
which implies that for $x\in (0,r)$
\begin{equation}
\label{eq:Ywithweight}
\begin{pmatrix} Y_{11,+}(x) \\ Y_{21,+}(x)\end{pmatrix} =\frac{(-1)^n}{\sqrt{\omega_{\gamma,n,\nu}^+(x)}}  e^{\frac{n\ell_\gamma}{2}\sigma_3} R(x) E_n(x) \begin{pmatrix}
J_\nu(2n\sqrt{f_\gamma(x)}) \\
2\pi i n \sqrt{f_\gamma(x)} J'_\nu(2n\sqrt{f_\gamma(x)})
\end{pmatrix}.
\end{equation}
Similarly, one obtains that
\begin{align*}
\begin{pmatrix}
-Y_{21}(z) \\ Y_{11}(z)
\end{pmatrix}^{T}
&= \begin{pmatrix} 0 & 1 \end{pmatrix} Y(z)^{-1} \\
&= z^{-\frac{\nu}{2}} e^{n(g_\gamma(z)-\frac{\ell_\gamma}{2}+\varphi_\gamma(z))} 
\begin{pmatrix} -1 & 1\end{pmatrix} e^{-\frac{\pi i \nu}{2}\sigma_3}
B(n^2 f_\gamma(z))^{-1} E_n(z)^{-1} R(z)^{-1} e^{-\frac{n\ell_\gamma}{2}\sigma_3}
\end{align*}
Inverting $B$ is straightforward, since it has determinant $1$. Using this and repeating the same steps as above, we obtain for $x\in (0,r)$ that
\begin{equation}
\label{eq:transposeYwithweight}
\begin{pmatrix}
-Y_{21,+}(x) \\ Y_{11,+}(x)
\end{pmatrix}^{T}
= \frac{(-1)^n}{\sqrt{\omega_{\gamma,n,\nu}^+(x)}}\begin{pmatrix}
2\pi i n \sqrt{f_\gamma(x)} J'_\nu(2n\sqrt{f_\gamma(x)}) \\
-J_\nu(2n\sqrt{f_\gamma(x)})
\end{pmatrix}^{T} E_n(x)^{-1} R(x)^{-1} e^{-\frac{n\ell_\gamma}{2}\sigma_3}.
\end{equation}
For any $x\in (0,\infty)$, we use the notation
\begin{equation}
\label{eq:defxn}
x_n=\frac{c_\gamma x}{\pi^2 n^2},
\end{equation}
for every $n\geq 1$, where $c_\gamma$ is defined as in \eqref{eq:defcgamma}. Note that for $n$ big enough, we have that $x_n\in (0,r)$ and that $x_n$ is precisely the quantity that appears in the scaling limit in \eqref{eq:resultRHprob}. 

Now let $(x,y) \in S$, where $S$ is a compact subset of the first quadrant. We turn to the asymptotic behaviour on $S$ as $n\rightarrow \infty$. From \eqref{eq:fgamma0} and the definition \eqref{eq:defcgamma} of $c_\gamma$ we infer that
\begin{align*}
2n\sqrt{f_\gamma(x_n)}&=\sqrt{x}+O\left(\frac{x^{3/2}}{n^2}\right).
\end{align*}
Since $x$ is bounded both from below and from above we thus have, using the mean value theorem, that uniformly
\begin{align} \label{eq:behavJnu}
J_\nu\left(2n\sqrt{f_\gamma(x_n)}\right) &= J_\nu(\sqrt{x})+O\left(\frac{1}{n^2}\right) \\ \label{eq:behavJnu'}
2n\sqrt{f_\gamma(x_n)} J'_\nu\left(2n\sqrt{f_\gamma(x_n)}\right) &= \sqrt x J'_\nu(\sqrt{x})+O\left(\frac{1}{n^2}\right).
\end{align}
By Cauchy's integral formula we have uniformly for $|\xi|<\frac{r}{3}$ that as $n\to\infty$
\begin{align} \label{eq:RderivCauchy}
R'(\xi) = \frac{1}{2\pi i} \oint_{|z|=\frac{r}{2}} \frac{R(z)-I}{(z-\xi)^2} dz = \mathcal O\left(\frac{1}{n}\right),
\end{align}
where in the last step we have used \eqref{eq:asymptoticR}. 
Then by \eqref{eq:defxn} and \eqref{eq:asymptoticR} again we have
\[R(y_n)^{-1} R(x_n) = I +R(y_n)^{-1}(R(x_n)-R(y_n))
= I + R(y_n)^{-1} \int_{y_n}^{x_n} R'(\xi) d\xi
= I+\mathcal O\left(\frac{x-y}{n^3}\right),\]
uniformly for $(x,y)\in S$. From the definition \eqref{eq:defEn} and RH-N4a it follows that $E_n(x_n)=\mathcal O(\sqrt{n})$ and $E_n(x_n)^{-1}=\mathcal O(\sqrt{n})$ as $n \rightarrow \infty$. Repeating the argument in \eqref{eq:RderivCauchy} for $E_n(z)$ instead of $R(z) - I$ then yields
\[E_n(y_n)^{-1} E_n(x_n) = I + \mathcal O\left(\frac{x-y}{n}\right),\]
as $n\to\infty$ uniformly for $(x,y)\in S$. We conclude that uniformly for $(x,y)\in S$
\begin{align} \label{eq:behavERRE}
E_n(y_n)^{-1} R(y_n)^{-1} R(x_n) E_n(x_n) &= I + \mathcal O\left(\frac{x-y}{n}\right) + \mathcal O\left(E_n(y_n)^{-1}\frac{x-y}{n^3} E_n(x_n)\right) \nonumber \\
&= I + \mathcal O\left(\frac{x-y}{n}\right), 
\end{align}
as $n\to\infty$. Here we have used again that $E_n(x_n)=\mathcal O(\sqrt{n})$ and $E_n(y_n)^{-1}=\mathcal O(\sqrt{n})$. We now assemble this all to prove Proposition \ref{prop:resultRHprob}. For this, we start with
\begin{align*}
\frac{c_\gamma}{\pi^2 n^2} K_n(x_n,y_n;\omega_{\gamma,n,\nu}^+) &= \frac{c_\gamma}{\pi^2 n^2} \frac{1}{2\pi i (x_n-y_n)} \sqrt{\omega_{\gamma,n,\nu}^+(x_n)\omega_{\gamma,n,\nu}^+(y_n)}  \begin{pmatrix} 0 & 1\end{pmatrix}
Y_+(y_n)^{-1} Y_+(x_n) 
\begin{pmatrix}
1\\ 0
\end{pmatrix} \\
&= \frac{1}{2\pi i (x-y)} \sqrt{\omega_{\gamma,n,\nu}^+(x_n)\omega_{\gamma,n,\nu}^+(y_n)}  \begin{pmatrix}
-Y_{21,+}(y_n) \\ Y_{11,+}(y_n)
\end{pmatrix}^{T} \begin{pmatrix} Y_{11,+}(x) \\ Y_{21,+}(x)\end{pmatrix},
\end{align*}
by \eqref{eq:kernelinY}. Then inserting \eqref{eq:Ywithweight} and \eqref{eq:transposeYwithweight}, we obtain
\begin{align*}
\frac{c_\gamma}{\pi^2 n^2} K_n(x_n,y_n;\omega_{\gamma,n,\nu}^+) = \frac{1}{2\pi i (x-y)}  & \begin{pmatrix}
2\pi i n \sqrt{f_\gamma(y_n)} J'_\nu(2n\sqrt{f_\gamma(y_n)}) \\
-J_\nu(2n\sqrt{f_\gamma(y_n)})
\end{pmatrix}^{T} E_n(y_n)^{-1} R(y_n)^{-1} \\
& \times R(x_n) E_n(x_n) \begin{pmatrix}
J_\nu(2n\sqrt{f_\gamma(x_n)}) \\
2\pi i n \sqrt{f_\gamma(x_n)} J'_\nu(2n\sqrt{f_\gamma(x_n)})
\end{pmatrix}.
\end{align*}
Inserting the asymptotic behaviours \eqref{eq:behavJnu}, \eqref{eq:behavJnu'} and \eqref{eq:behavERRE}, we obtain, as $n\to \infty$, that 
\begin{align*}
\frac{c_\gamma}{\pi^2 n^2} K_n(x_n,y_n;\omega_{\gamma,n,\nu}^+) &= \frac{1}{2\pi i (x-y)}\begin{pmatrix}
\pi i \sqrt{y} J'_\nu(\sqrt{y}) \\
-J_\nu(\sqrt{y})
\end{pmatrix}^{T} \begin{pmatrix}
J_\nu(\sqrt{x}) \\
\pi i \sqrt{x} J'_\nu(\sqrt{x})
\end{pmatrix} + \mathcal O\left(\frac{1}{n}\right) \\
&= \frac{J_\nu(\sqrt{x})\sqrt{y}J'_\nu(\sqrt{y})-J_\nu(\sqrt{y})\sqrt{x}J'_\nu(\sqrt{x})}{2(x-y)} + \mathcal O\left(\frac{1}{n}\right),
\end{align*}
uniformly for $(x,y)\in S$. This concludes the proof. 
\end{proof}